\documentclass[11pt,reqno,letterpaper]{amsart}
\usepackage{fix-cm}
\usepackage{orcidlink}
\usepackage{amsfonts}
\usepackage{graphicx}
\usepackage{verbatim}
\usepackage{amscd}
\usepackage{amsmath, physics}
\usepackage{amssymb}
\usepackage{latexsym}
\usepackage{hyperref}
\usepackage[all]{xy}
\usepackage{color}
\usepackage{mathrsfs}
\newtheorem{theorem}{\bf Theorem}
\newtheorem{proposition}{\bf Proposition}

\newtheorem{lemma}{\bf Lemma}
\theoremstyle{definition}
\newtheorem{remark}{\bf Remark}
\newtheorem{definition}{\bf Definition}
\newtheorem{example}{\bf Example}

\usepackage{empheq,etoolbox}
\numberwithin{equation}{section}
\patchcmd{\subequations}
  {\theparentequation\alph{equation}}
  {\theparentequation.\alph{equation}}
  {}{}
  
\begin{document}
\title[MCF into an ambient Riemannian manifold evolving by R--B flow]{Mean curvature flow into an ambient Riemannian manifold evolving by Ricci--Bourguignon flow}

\author[José N.V. Gomes]{José Gomes$^1$\orcidlink{0000-0001-5678-4789}}
\author[Matheus Hudson]{Matheus Hudson$^2$\orcidlink{0000-0001-9786-1234}}
\author[Raquel Lemos]{Raquel Lemos$^3$}

\address{$^{1}$Department of Mathematics, Universidade Federal de São Carlos, Rod. Washington Luís, Km 235, 13565-905, São Carlos, São Paulo, Brazil.}
\address{$^{2}$Department of Mathematics, Universidade Federal do Amazonas, Av. General Rodrigo Octavio Jordão Ramos, 1200, 69080-900, Manaus, Amazonas, Brazil.}
\address{$^3$Institute of Mathematics and Computer Science, Universidade de São Paulo, Avenida Trabalhador Sancarlense, 400, 13566-590, São Carlos, São Paulo, Brazil.}

\email{$^{1}$jnvgomes@ufscar.br}
\email{$^{2}$matheushudson@ufam.edu.br}
\email{$^{3}$raquellemos@usp.br}

\urladdr{$^{1}$https://www.ufscar.br}
\urladdr{$^{2}$https://ufam.edu.br}
\urladdr{$^{3}$https://www.icmc.usp.br}


\begin{abstract}
We study the mean curvature flow in Riemannian manifolds evolving by the Ricci--Bourguignon flow. Using DeTurck's trick, we establish short-time existence and uniqueness for the mean curvature flow in this time-dependent ambient geometry. In the case where the ambient metric evolves by a self-similar solution to the Ricci--Bourguignon flow, we characterize the associated mean curvature solitons and derive a Huisken monotonicity-type formula adapted to the Ricci--Bourguignon setting. As a consequence, we show that, in compact ambient manifolds, type-I singularities of the mean curvature flow in a gradient shrinking Ricci--Bourguignon soliton background have, after suitable parabolic rescaling, weighted minimal hypersurfaces as their limiting models, while in the noncompact setting, we establish the corresponding result for a gradient shrinking Schouten soliton background. This provides a geometric framework for understanding how the ambient scalar curvature influences the formation and structure of singularity models for the mean curvature flow. We also show how to construct explicit examples of mean curvature flow in a gradient Ricci--Bourguignon soliton background.
\end{abstract}
\maketitle

\section{\bf Introduction}
The mean curvature flow has found applications in various fields, including geometric analysis, geometric measure theory, and partial differential equations. Historically, Brakke~\cite{Brakke} pioneered the study of how surfaces deform due to their mean curvature within the framework of geometric measure theory, as noted by Huisken~\cite{Huisken}. One of the motivations for studying such flows is that these moving surfaces can be used to describe the behavior of grain boundaries in pure metals during the heating process. Huisken approached this problem from a more classical perspective by studying convex surfaces evolving under mean curvature flow in Euclidean space. He showed that these surfaces become asymptotically spherical and eventually shrink to a single point in finite time. 

In 1986, Gage and Hamilton~\cite{gage1986heat} introduced a generalization of the mean curvature flow to immersions into ambient Riemannian manifolds. Since then, the mean curvature flow has been a constant object of investigation and has experienced great development in the last few decades. In the setting where an ambient Riemannian manifold evolves by Ricci flow, it is worth mentioning the works of Lott~\cite{lott2012mean} when studying the mean curvature flow in a Ricci flow background, and of Magni, Mantegazza and Tsatis~\cite{magni2013flow} when studying the mean curvature flow of a submanifold inside a moving Riemannian manifold. 

Subsequently, Yamamoto~\cite{yamamoto2020meancurvature} proved that the mean curvature flow in a gradient shrinking Ricci soliton background converges, after rescaling, to a self-shrinker when it develops a singularity of type-I. Very recently, the first two authors, together with Yamamoto~\cite{gomes2024anote}, applied a Huisken monotonicity-type formula to the mean curvature flow in an ambient space whose Riemannian metric evolves by a shrinking self-similar solution of an extended Ricci flow. 

Here, we find a Huisken monotonicity-type formula for the mean curvature flow in an ambient smooth manifold whose Riemannian metric evolves by a self-similar solution to the Ricci--Bourguignon flow. We also apply such a formula in the shrinking case.

The Ricci--Bourguignon flow, introduced by Bourguignon~\cite{Bour}, is a smooth one-parameter family of Riemannian metrics $\{g(t)\}_{t\in[0,T)}$ on an $(n\geq 3)$-dimensional smooth manifold $M$ evolving according to the equation
\begin{equation} \label{eq-RB}
\frac{\partial}{\partial t} g(t)=-2\big(\operatorname{Ric}_{g(t)}-\rho R g(t)\big) \,\,\, \forall t\in[0,T),
\end{equation}
where $\operatorname{Ric}_{g(t)}$ is the Ricci tensor of the Riemannian manifold $(M,g(t))$, $R$ its scalar curvature and $\rho$ is a real constant. Note that $\rho=0$ corresponds to the Ricci flow introduced by Hamilton~\cite{Hamilton92}, while $\rho=1/2(n-1)$, up to a constant rescaling of time, corresponds to the Schouten flow.

Catino et al.~\cite{Catino} developed a parabolic theory for the Ricci--Bourguignon flow, proving short-time existence and curvature estimates. More precisely, when $\rho < {1}/{2(n-1)}$, each initial compact Riemannian manifold $(M, g_0)$ without boundary has a unique smooth solution $g(t)$ satisfying the Ricci--Bourguignon flow equation, at least for a short time interval. Most recently, Catino and Mantegazza~\cite{catino2026short} proved short-time existence and uniqueness for the Schouten flow.

A gradient Ricci--Bourguignon soliton is a self-similar solution $\overline{g}(t)$ of the Ricci--Bourguignon flow given by
\begin{equation*}
\overline{g}(t)  = \sigma(t)\psi_t^* g,
\end{equation*}
for some initial Riemannian metric $g$, where $\psi_t$ is a  smooth one-parameter family of diffeomorphisms of $M$ generated from the flow of $\nabla_g f/\sigma(t)$, for some $f\in C^\infty(M)$, and  $\sigma(t)$ is a smooth positive function on $t$. By setting $\overline{f}(t)=\psi_t^*f$, it follows from~\eqref{eq-RB} that
\begin{equation}\label{mod_grad_Ricci_soliton}
\operatorname{Ric}_{\overline{g}} + \operatorname{Hess}_{\overline{g}}\,\overline{f} - \rho R_{\overline{g}}\overline{g} = \dfrac{c}{2(T - t)}\overline{g},
\end{equation}
where $c= 0$ in the steady case, for $t \in  \mathbb{R};$ $c = 1$ in the shrinking case, for $t \in (-\infty, T)$;  and $c = 1$ in the expanding case, for $t \in (T,+\infty).$ Moreover,
\begin{align}\label{self-solution}
\dfrac{\partial }{\partial t} \overline{f} = \left|\nabla_{\overline{g}} \overline{f}\right|_{\overline{g}}^2\,.
\end{align}
The function $\overline{f}$ is the so-called \emph{potential function}. 

In this article, we consider the mean curvature flow (MCF) of an orientable compact $(n-1)$-dimensional smooth manifold $\Sigma$ without boundary in an $n$-dimensional smooth ambient manifold $M$ equipped with a complete Riemannian metric that evolves under a Ricci--Bourguignon flow. More precisely, let $g(t)$ be a Ricci--Bourguignon flow on $M$, and let $\{F(\cdot,t); \; t\in [0,T)\}$ be a smooth one-parameter family of immersions of $\Sigma$ in $M$.  For each $t\in [0,T),$ set  $F_t = F(\cdot, t)$ and $\Sigma_t$ for the hypersurface $F_t(\Sigma)$ of $\big(M,g(t)\big),$ which can be written as $\Sigma_t:=\big(\Sigma, F^*_t g(t)\big)$. We say that the family $\mathscr F :=\{\Sigma_t;\; t\in [0,T)\}$ is the \emph{MCF in a Ricci--Bourguignon flow background} if it evolves under MCF
\begin{align}\label{mean:curvature:RB-flow}
\left\{
\begin{array}{rcl}
\dfrac{\partial}{\partial t}F(p,t) &=&  H(p,t) \nu(p,t),\\ [1ex]
F(p, 0) &=& F_0(p),
\end{array}
\right.
\end{align}
where $H(p,t)$ and $\nu(p,t)$ are the mean curvature and the unit normal of $\Sigma_t$ at the point $p\in\Sigma$, respectively. In the particular case where $\overline{g}(t)$ is a gradient Ricci--Bourguignon soliton on $M$, with potential function $\overline{f}$, we say that $\mathscr F$ is the \emph{MCF in a gradient Ricci--Bourguignon soliton background}. In the special case  $\rho=1/2(n-1)$, we say that $\mathscr F$ is the \emph{MCF in a gradient Schouten soliton background}. A hypersurface $\Sigma_t\in \mathscr F$ is called a \emph{mean curvature soliton}, if
\begin{align*}
H(p,t)+\nu(p,t)\overline{f}=0 \quad \forall p\in\Sigma.
\end{align*}
Here, $\nu(\,\cdot\,,t)$ must be the inward unit normal vector field on $\Sigma_t.$

In Proposition~\ref{existence:MCF-Ricci Bourguignon:flow} we show short-time existence and uniqueness of the MCF in the Ricci--Bourguignon flow background. 

In Section~\ref{Sec-CMCS}, we show how to construct a family $\mathscr{G}$ of mean curvature solitons. We prove that a sufficient condition for $\mathscr{G}$ to be a family of mean curvature solitons is that the hypersurface $\Sigma$ be $f$-minimal. Moreover, we show that any family $\mathscr{F}$ of mean curvature solitons is given by $\mathscr{G}$, up to reparametrization. This is the content of our first main theorem.

\begin{theorem}\label{NazasCharacterization}
Let $M$ be an $(n\geq3)$-dimensional smooth manifold, and let $\overline{g}(t)$ be a gradient Ricci--Bourguignon soliton on $M$ for some initial value $g$. If $\Sigma$ is an $f$-minimal hypersurface of $(M,g)$, then $\mathscr G$, defined as in Section~\ref{Sec-CMCS}, is a family of mean curvature solitons. Moreover, any family $\mathscr F$ of mean curvature solitons is given by $\mathscr G$ up to reparametrization.
\end{theorem}

Theorem~\ref{NazasCharacterization} provides an essential tool for constructing explicit examples of MCF in a gradient Ricci--Bourguignon soliton background. In Section~\ref{example:section}, we establish Proposition~\ref{Feitosa-FreitasFilho-Gomes}, which provides another useful tool for this construction.

Combining the evolution of the induced area element under the Ricci--Bourguignon flow background (see Lemma~\ref{Huisken:lem}) with Hamilton's approach~\cite{Hamilton99}, an estimate due to Catino et al.~\cite{Catino}, and the Gauss--Bonnet theorem, we obtain the following finite-time singularity result for the MCF of two-spheres in a Ricci--Bourguignon flow background.

\begin{theorem}\label{Cor-FTS}
Let $M$ be a three-dimensional compact smooth manifold, and let $\{g(t)\}_{t\in[0,T)}$ be a maximal solution of the Ricci--Bourguignon flow, with $\rho \leq1/4$. If $\{F_t\}_{t\in[0, S)}$, $0<S\leq T$, is a maximal solution of the MCF of the two-sphere in a Ricci--Bourguignon flow background on $M$, then $S<\infty$. In particular, the ambient flow and the MCF cannot both remain smooth for all positive times.
\end{theorem}

The next result is a Huisken monotonicity-type formula for the MCF in a gradient Ricci--Bourguignon soliton background.

\begin{theorem}\label{Huisken_monotonicity}
Let $M$ be an $(n\geq3)$-dimensional smooth manifold, and let $\overline{g}(t)$ be a gradient soliton to the Ricci--Bourguignon flow on $M$ with the potential function~$\overline{f}$. Assume that $\mathscr{F}$ is the MCF in a gradient Ricci--Bourguignon soliton background, and denote by $dA_{\overline{g}}$ the $(n-1)$-dimensional Riemannian measure on $\Sigma_t$ induced by $\overline g(t).$ If $\Phi(t)$ is the smooth function given by:
\begin{itemize}
\item[(i)]  $\mathbb{R}\ni t\mapsto\int_\Sigma e^{-\overline{f}}dA_{\overline{g}}$ in the steady case;
\smallskip
\item[(ii)] $(-\infty,T)\ni t\mapsto [4\pi(T-t)]^{\frac{-(n-1)}{2}}\int_\Sigma e^{-\overline{f}}dA_{\overline{g}}$ in the shrinking case;
\smallskip
\item[(iii)] $(T,\infty)\ni t\mapsto [4\pi(t-T)]^{\frac{-(n-1)}{2}}\int_\Sigma e^{-\overline{f}}dA_{\overline{g}}$ in the expanding case.
\smallskip
\end{itemize}
Then, $\Phi(t)$ is nonincreasing. Moreover, it is constant if and only if $\mathscr F$ is a family of mean curvature solitons.
\end{theorem}

The most appropriate setting to apply Theorem~\ref{Huisken_monotonicity} is to consider the normalized family $\widetilde{\mathscr{F}}$ of $\mathscr{F}$ in a gradient shrinking Ricci--Bourguignon soliton background as follows. 

Let $\mathscr F :=\{\Sigma_t; \; t\in [0,T)\}$ be the MCF in a gradient shrinking Ricci--Bourguignon soliton background with $\overline{g}(t)=(T-t)\psi_t^{*}g$, for $t\in (-\infty,T)$. Setting $s=-\log (T- t)$ and $\widetilde{\Sigma}_s = (\Sigma,\widetilde{F}_{s}^{*}g)$, where $\widetilde{F}_{s}=\psi_t\circ F_t$, we obtain
\begin{equation*}
\frac{\partial \widetilde F_s}{\partial s} = \big(\nabla_g f + H_g\big)(\widetilde F_s).
\end{equation*}
We call $\widetilde{\mathscr F}:=\{\widetilde \Sigma_s;\; s\in [-\log T,\infty)\}$ the \emph{normalized MCF} in $(M,g)$.

The MCF in a gradient shrinking Ricci--Bourguignon soliton background develops a \emph{singularity of type-I} when there exists a constant $C > 1$ such that 
\begin{equation}\label{singularity:type:I}
\max_{p \in \Sigma} |\mathcal A(p, t)| \leq\frac{C}{\sqrt{T - t}}, 
\end{equation}
where $\mathcal A(\cdot, t)$ stands for the second fundamental form of $\Sigma_t$.

We apply the Huisken monotonicity-type formula for the MCF in a gradient shrinking Ricci--Bourguignon soliton background when $\rho\neq 0$. In this setting, it is natural to consider a function $\phi\in C^\infty(M)$ satisfying $d\phi=Rdf$, as follows from the Hamilton-type equation (see Section~\ref{sect:apllication}). Using this function, we establish some preliminary results that enable us to prove our next main theorem. More precisely, we show that, in compact ambient manifolds, type-I singularities of the MCF in a gradient shrinking Ricci--Bourguignon soliton background have, after suitable parabolic rescaling, weighted minimal hypersurfaces as their limiting models, while in the noncompact setting, we establish the corresponding result for the MCF in a gradient shrinking Schouten soliton background. The case $\rho=0$, corresponding to the Ricci flow, was previously treated by Yamamoto~\cite{yamamoto2020meancurvature}, without requiring the assumption on $\phi$.

\begin{theorem}\label{application:Huisken}
Let $(M,g)$ be an $(n\geq 3)$-dimensional compact Riemannian manifold, and let $\mathscr F$ be the MCF in a gradient shrinking Ricci--Bourguignon soliton background which develops a singularity of type-I. Consider the normalized MCF $\widetilde{\mathscr F}$ in $(M,g)$ such that $\phi\circ \widetilde{F}_{s},$ with $\rho \leq 1/2(n-1)$ and $\rho\neq0,$ remains uniformly bounded on $[-\log T,\infty).$  Let $\left\{s_{j}\right\}_{j=1}^{\infty}$ be a divergent increasing  sequence, and let $\left\{p_{j}\right\}_{j=1}^{\infty}$ be a sequence of points in $\Sigma.$ Then, there exist subsequences $\{s_{j_{k}}\}$ and $\{p_{j_{k}}\}$ such that the family of immersion maps $\widetilde F_{s_{j_k}}: \Sigma \to (M,g)$ from pointed manifolds $\big(\Sigma,  p_{j_{k}}\big)$ converges to an immersion map $F_\infty: \Sigma_\infty \to (M,g)$ from an $(n-1)$-dimensional complete pointed Riemannian manifold $\left(\Sigma_{\infty}, F_\infty^*g, p_{\infty}\right)$ in the $C^\infty$ Cheeger--Gromov sense. Furthermore, $(\Sigma_{\infty}, F_\infty^* g)$ is an $f_\infty-$minimal hypersurface of  $(M, g),$ where $f_\infty = f \circ F_\infty$. 
\end{theorem}

In the noncompact setting, the Schouten flow provides a suitable geometric framework for establishing convergence in the Cheeger--Gromov sense. Under appropriate geometric and analytic assumptions, Theorem~\ref{Thmnoncommpact} establishes the corresponding Cheeger--Gromov convergence of the MCF in a gradient shrinking Schouten soliton background on a noncompact manifold.

\section{\bf Short-time existence and uniqueness for the MCF in a Ricci--Bourguignon flow background}\label{SecMCFERFB}

The following proposition establishes short-time existence and uniqueness for the MCF in a Ricci--Bourguignon flow background. Our approach follows the DeTurck method employed by Han--Li--Zhao~\cite[Theorem~2.1]{han2018mean} to establish short-time existence for the MCF along a K\"ahler--Ricci flow. Since the argument relies only on the smooth dependence of the ambient metric on space and time, it extends directly to any prescribed smooth family of Riemannian metrics. We provide the details of the proof, including the uniqueness argument, which was not addressed in~\cite{han2018mean}.

\begin{proposition}\label{existence:MCF-Ricci Bourguignon:flow}
Let $\{g(t)\}_{t\in[0,T)}$ be a solution to the Ricci--Bourguignon flow on an $(n\geq3)$-dimensional smooth manifold $M$, and let $F:\Sigma\to M$ be an immersion of an orientable compact $(n-1)$-dimensional smooth manifold $\Sigma$ without boundary. Then there exists a unique smooth one-parameter family of immersions $F_t:\Sigma\to M$, with $F_0=F$, evolving according to~\eqref{mean:curvature:RB-flow} with respect to the induced metric $F_t^*g(t)$ on $\Sigma$, for some short-time interval.
\end{proposition}

\begin{proof}
Choose local coordinates $(x_1,\dots,x_{n-1})$ on $\Sigma$, and let $F^\alpha(x,t)$ be the coordinate components of $F(x,t)$ in a local coordinate system of $M\times [0,T)$, with $F(x,0)=F(x)$, for some $F(x,t)$ to be determined. Consider $\widehat\Gamma^k_{ij}(t)$ and $\Gamma^\alpha_{\lambda\beta}(t)$ the Christoﬀel symbols of $F^*_tg(t)$ and of $g(t),$ respectively, where $F_t(x):=F(x,t)$. In what follows, the subscript \(``\widehat{\phantom{X}}"\) indicates that the corresponding geometric quantities are computed with respect to the metric $F_t^*g(t)$. Thus, we have 
\begin{equation*}
dF_t\left(\frac{\partial}{\partial x^i}\right) = \frac{\partial F_t^\alpha}{\partial x^ i}\frac{\partial}{\partial x^\alpha} \quad
\mbox{and} \quad dF_{t}\left(\widehat{\nabla}_{\frac{\partial}{\partial x^i}}\frac{\partial}{\partial x^j}\right) = \widehat\Gamma^k_{ij}(t)\frac{\partial F_{t}^\alpha}{\partial x^k}\frac{\partial}{\partial x^\alpha}.
\end{equation*}
Furthermore, $\frac{\partial}{\partial x^i}\left(\frac{\partial F_{t}^\alpha}{\partial x^ j}\circ F_{t}\right)=dF_{t}\left(\frac{\partial}{\partial x^ i}\right)\left(\frac{\partial F_{t}^\alpha}{\partial x^ j}\right)$ (see, e.g., \cite{Lee} p.~61). Hence,
\begin{align*}
\nabla_{dF_t\left(\frac{\partial}{\partial x^ i}\right)}dF_t\left(\frac{\partial}{\partial x^j}\right) = \left(\frac{\partial^{2}F_{t}^\alpha}{\partial x^i\partial x^j} + \Gamma^\alpha_{\lambda\beta}(t)\frac{\partial F_{t}^\lambda}{\partial x^i}\frac{\partial F_{t}^\beta}{\partial x^j}\right) \frac{\partial}{\partial x^\alpha}
\end{align*}
and the Gauss equation
\begin{equation*}
\nabla_{dF_t\left(\frac{\partial}{\partial x^i}\right)}dF_t\left(\frac{\partial}{\partial x^j}\right) = dF_{t}\left(\widehat{\nabla}_{\frac{\partial}{\partial x^i}}\frac{\partial}{\partial x^j}\right) + \mathcal{A}_{ij}\nu
\end{equation*}
reduces to
\begin{align*}
\left(\frac{\partial^2 F_t^\alpha}{\partial x^i \partial x^j} + \Gamma^{\alpha}_{\lambda \beta}(t) \frac{\partial F_t^{\lambda}}{\partial x^i}\frac{\partial F_t^{\beta}}{\partial x^j}- \widehat\Gamma^{k}_{ij}(t)\frac{\partial F_t^{\alpha}}{\partial x^{k}}\right)\frac{\partial}{\partial x^\alpha} = \mathcal{A}_{ij}\nu.
\end{align*}
Moreover, it is well-known that
\begin{align*}
\widehat\Delta F_t = \big(F_t^*g(t)\big)^{ij}\widehat\nabla_i\widehat\nabla_j F_t = \big(F_t^* g(t)\big)^{ij}\left( \frac{\partial^2 F_t}{\partial x^i \partial x^j} - \widehat\Gamma^k_{ij}(t)\frac{\partial F_t}{\partial x^k}
\right).
\end{align*}
Consequently, taking the trace of the preceding Gauss equation, the family $F_t$ evolves under MCF \eqref{mean:curvature:RB-flow} if
\begin{equation}\label{nonparabolicequation}
\frac{\partial F_t^\alpha}{\partial t}=\widehat\Delta F_t^\alpha+\big(F_t^*g(t)\big)^{ij}\,\Gamma^\alpha_{\lambda\beta}(t)\frac{\partial F_t^\lambda}{\partial x^i}\frac{\partial F_t^\beta}{\partial x^j}.
\end{equation}
The problem is that \eqref{nonparabolicequation} is not a strictly parabolic equation. Indeed, by a straightforward calculation, we can deduce that the principal symbol is
\begin{equation*}\big(F^*_tg(t)\big)^{ij}\zeta_i\zeta_j\left[v^\alpha- \big(F^*_tg(t)\big)^{k\ell} v^\lambda\frac{\partial F_t^\beta}{\partial x^\ell }g_{\lambda\beta}\frac{\partial F_t^\alpha}{\partial x^k }\right],
\end{equation*}
for all cotangent vector $\zeta$, which clearly shows the existence of zeros in this symbol. Fortunately, we can use DeTurck's trick to modify the flow~\eqref{nonparabolicequation} through a diffeomorphism of the parameter space of $\Sigma$. For this, we consider the following equation
\begin{equation}\label{evolution+of+tilde}
\frac{\partial \widetilde F_t^\alpha}{\partial t}
\!=\!\big(\widetilde F^*_t g(t)\big)^{ij}\!\!\left(\frac{\partial^2\widetilde F_t^\alpha}{\partial y^i\partial y^j}\! -\! \widetilde\Gamma_{ij}^k(t) \frac{\partial \widetilde F_t^\alpha}{\partial y^k}+{\Gamma}^\alpha_{\lambda\beta}(t)
\frac{\partial \widetilde F_t^\lambda}{\partial y^i}
\frac{\partial \widetilde F_t^\beta}{\partial y^j}  \right) \!+\!
V^k \frac{\partial\widetilde F_t^\alpha}{\partial y^k}
\end{equation}
where $\widetilde F(y(x, t), t)=F(x, t)$, and the vector field $V^k$ will be chosen to make the equation strictly parabolic. Furthermore, we are using the subscript $``\widetilde{\phantom{X}}"$ to denote the geometric objects computed with respect to the metric $\widetilde F_t^*g(t)$. 

We note that if $\widetilde F(y(x, t), t)$ is a solution of~\eqref{evolution+of+tilde}, then 
\begin{align*}
\frac{\partial F_t^\alpha}{\partial t} &= \frac{\partial \widetilde F_t^\alpha}{\partial t}  + \frac{\partial \widetilde F_t^\alpha}{\partial y^k} \frac{\partial y^k}{\partial t} \\
&= \widetilde\Delta \widetilde F_t^\alpha + \big(\widetilde F^*_t g(t)\big)^{ij}{\Gamma}^\alpha_{\lambda\beta}(t)
\frac{\partial \widetilde F_t^\lambda}{\partial y^i}
\frac{\partial \widetilde F_t^\beta}{\partial y^j} 
+ \left( V^k + \frac{\partial y^k}{\partial t} \right) \frac{\partial\widetilde F_t^\alpha}{\partial y^k} \\
&=\widehat\Delta  F_t^\alpha + \big(F^*_tg(t)\big)^{ij}{\Gamma}^\alpha_{\lambda\beta}(t)
\frac{\partial F_t^\lambda}{\partial x^i}
\frac{\partial F_t^\beta}{\partial x^j}.
\end{align*}
Thus, we recover \eqref{nonparabolicequation}, since we can consider
\begin{equation*}
\begin{cases}
\dfrac{\partial y^k}{\partial t} = -V^k(y, t), \\
y_0 = \operatorname{Id}_\Sigma,
\end{cases} 
\end{equation*}
where 
\begin{equation*}
V^k = \big(\widetilde F^*_t g(t)\big)^{ij} \left( \widetilde\Gamma^k_{ij}(t) - \widetilde\Gamma^k_{ij}(0) \right).
\end{equation*}
Furthermore, equation~\eqref{evolution+of+tilde} reduces to
\begin{align*}
\frac{\partial \widetilde F_t^\alpha}{\partial t} \!=& \big(\widetilde F^*_t g(t)\big)^{ij}\!\! \left(\!\! \frac{\partial^2 \widetilde F_t^\alpha}{\partial y^i \partial y^j}\!-\! \widetilde\Gamma^{k}_{ij}(t)\frac{\partial\widetilde F_t^\alpha}{\partial y^k} \!\!\right)\!\!+\!\big(\widetilde F^*_t g(t)\big)^{ij}\Gamma^{\alpha}_{\lambda\beta}(t) \frac{\partial \widetilde F_t^\lambda}{\partial y^i}
\frac{\partial \widetilde F_t^\beta}{\partial y^j} \\
&+\big(\widetilde F^*_t g(t)\big)^{ij} \left( \widetilde\Gamma^k_{ij}(t) -\widetilde \Gamma^k_{ij}(0) \right)\!\!\frac{\partial\widetilde F_t^\alpha}{\partial y^k}\\
=& (\widetilde F^*_tg(t))^{ij} \left( \frac{\partial^2 \widetilde F_t^\alpha}{\partial y^i \partial y^j} - \widetilde\Gamma^k_{ij}(0) \frac{\partial \widetilde F_t^\alpha}{\partial y^k} \right)+\big(\widetilde F^*_t g(t)\big)^{ij}\Gamma^{\alpha}_{\lambda\beta}(t) \frac{\partial \widetilde F_t^\lambda}{\partial y^i}
\frac{\partial \widetilde F_t^\beta}{\partial y^j} ,
\end{align*}
whose principal symbol is $\big(\widetilde F^*_t g(t)\big)^{ij} \zeta_i\zeta_j \operatorname{Id}$. Hence, this modified evolution equation is strictly parabolic. Applying the theory of quasilinear parabolic systems (see, e.g., Aubin~\cite{Aubin}), we obtain a unique smooth short-time solution $\widetilde F_t$ with $\widetilde F_0=F$. Solving the preceding ODE for $y_t$ and setting $F_t=\widetilde F_t\circ y_t$, we obtain a smooth short-time solution of \eqref{mean:curvature:RB-flow} with initial condition $F_0=F$.

To prove uniqueness, let $G_t$ be another smooth solution of \eqref{mean:curvature:RB-flow} with $G_0=F$, and fix the reference
metric $h_0=F^*g(0)$. Consider the harmonic map heat flow
\begin{equation*}
\partial_t z_t=\tau_{G_t^*g(t),h_0}(z_t),\qquad z_0=\operatorname{Id}_{\Sigma},   
\end{equation*}
where $\tau_{G_t^*g(t),h_0}$ denotes the tension field for maps from $(\Sigma,G_t^*g(t))$ to $(\Sigma,h_0)$. Since $\Sigma$ is compact without boundary and $G_t^*g(t)$ depends smoothly on time, this strictly parabolic system admits a smooth solution for a sufficiently short time (see  Huang--Tam~\cite{HuangTam}). 

Now, set $\widetilde G_t=G_t\circ z_t^{-1}$ and define
\begin{equation*}
V_G^k = \bigl(\widetilde G_t^*g(t)\bigr)^{ij} \left(\Gamma(\widetilde G_t^*g(t))^k_{ij}-\widetilde\Gamma^k_{ij}(0)\right),    
\end{equation*}
where $\widetilde\Gamma^k_{ij}(0)$ are the Christoffel symbols of $h_0$. Using
\begin{equation*}
G_t^*g(t)=z_t^*\bigl(\widetilde G_t^*g(t)\bigr),    
\end{equation*}
we claim that the following identity holds:
\begin{align}\label{VG}
(\partial_tz_t)\circ z_t^{-1}=\tau_{\widetilde G_t^*g(t),h_0}(\operatorname{Id}_{\Sigma})=-V_G.    
\end{align}
Indeed, fix $t$ and write $\kappa=\widetilde G_t^*g(t).$ For a map
\begin{equation*}
z:(\Sigma,\kappa)\longrightarrow(\Sigma,h_0),
\end{equation*}
the tension field has components
\begin{equation*}
\tau_{\kappa,h_0}(z)^\alpha= \kappa^{ij}\left(\partial_i\partial_jz^\alpha-\Gamma(\kappa)^\ell_{ij}\partial_\ell z^\alpha+\Gamma(h_0)^\alpha_{\beta\gamma}(z)
\,\partial_i z^\beta\,\partial_j z^\gamma\right).    
\end{equation*}
Here $i,j,\ell,\alpha, \beta, \gamma \in\{1,\ldots,n-1\}.$ For \(z=\operatorname{Id}_\Sigma\), choose the same coordinates on the domain and target. Then
\begin{equation*}
z^\alpha(x)=x^\alpha,\qquad\partial_i z^\alpha=\delta_i^\alpha,\qquad\partial_i\partial_jz^\alpha=0.
\end{equation*}
Consequently,
\begin{align*}
\tau_{\kappa,h_0}(\operatorname{Id}_\Sigma)^\alpha
&=\kappa^{ij}\left(-\Gamma(\kappa)^\alpha_{ij}+\Gamma(h_0)^\alpha_{ij}\right)\\
&=-\bigl(\widetilde G_t^*g(t)\bigr)^{ij}
\left(\Gamma\big(\widetilde G_t^*g(t)\big)^\alpha_{ij}-\widetilde\Gamma^\alpha_{ij}(0)\right)\\
&=-V_G^\alpha,
\end{align*}
because $\widetilde\Gamma^\alpha_{ij}(0)=\Gamma(h_0)^\alpha_{ij}.$ It is enough to prove our claim. 

Next, note that the identity map need not have zero tension: its domain and target carry different metrics, hence different Levi–Civita connections. Writing $\vec H_{G_t,g(t)}$ for the mean curvature vector of $G_t$ in $(M,g(t))$, using~\eqref{VG}, and~\eqref{mean:curvature:RB-flow}, and the derivative
\begin{align*}
\partial_t G_t = (\partial_t\widetilde G_t) \circ z_t + d\widetilde G_t|_{z_t}(\partial_tz_t)
\end{align*}
of $G_t=\widetilde G_t\circ z_t,$ we obtain
\begin{align*}
\partial_t\widetilde G_t &= (\partial_tG_t)\circ z_t^{-1} -d\widetilde G_t \bigl((\partial_tz_t)\circ z_t^{-1}\bigr)\\
&=\vec  H_{\widetilde G_t,g(t)} + d\widetilde G_t(V_G).
\end{align*}
Therefore, $\widetilde G_t$ satisfies the same strictly parabolic modified system as $\widetilde F_t$, with the same
reference metric $h_0$ and the same initial condition $F$. Uniqueness for that system implies
\begin{equation*}
\widetilde G_t=\widetilde F_t.
\end{equation*}
Hence, $V_G=V$, and $z_t$ satisfies
\begin{equation*}
\partial_tz_t=-V(\,\cdot\,,t)\circ z_t,\qquad z_0=\operatorname{Id}_{\Sigma}.
\end{equation*}
By uniqueness for this ODE, $z_t=y_t$. Hence
\begin{equation*}
G_t=\widetilde G_t\circ z_t =\widetilde F_t\circ y_t =F_t.
\end{equation*}
This proves short-time uniqueness and completes the proof.
\end{proof}

\begin{remark}
The reason for restricting to $n\geq3$ is that, for $n=2$, identity \eqref{mod_grad_Ricci_soliton} reduces to the Yamabe soliton equation, while \eqref{eq-RB} reduces to a constant rescaling of the Yamabe flow whenever $\rho\neq1/2.$ In this case, the Ricci--Bourguignon flow does not exhibit any essentially new behavior compared with the Yamabe flow. For related results on the Yamabe flow, see Chow~\cite{Chow1992}, Daskalopoulos and Sesum~\cite{daskalo}, and Ye~\cite{Ye1994}. Nevertheless, we observe that Proposition~\ref{existence:MCF-Ricci Bourguignon:flow} remains valid in dimension $n=2$. \end{remark}

\section{\bf Proofs of Theorems~\ref{NazasCharacterization},~\ref{Cor-FTS} and~\ref{Huisken_monotonicity}}\label{Sec-CMCS}
In this section, we first show how to construct a family of mean curvature solitons and characterize such a family. For it, let $M$ be an $(n\geq3)$-dimensional smooth manifold, and let $\overline g(t)=\sigma(t)\psi_t^*g$ be a gradient Ricci--Bourguignon soliton on $M$ for some initial metric $g$ and potential function $\overline f=\psi^*_tf$, where $\{\psi_t\}$ is the smooth one-parameter family of diffeomorphisms of $M$ generated by $\nabla_g f/\sigma(t)$ with
\begin{enumerate}
\item[(i)] $\sigma(t) = \kappa(T-t)$ and $\psi_{T-\kappa} = \operatorname{Id}$, where $\kappa = 1$ in the shrinking case, for $t \in (-\infty, T$);
\smallskip
\item[(ii)]  $\sigma(t) = \kappa(T-t)$ and $\psi_{T-\kappa} = \operatorname{Id}$, where $\kappa = -1$ in the expanding case, for $t \in (T,+\infty)$;
\smallskip
\item[(iii)] $\sigma(t) = 1$ and $\psi_{0} = \operatorname{Id}$ in the steady case, for $t \in  \mathbb{R}$.
\end{enumerate}

Given an orientable compact $(n-1)$-dimensional smooth manifold $\Sigma$ without boundary in $M.$ We set $s(t)=T-\kappa/\sigma(t)$, $F_t = F(\,\cdot\,, t)$ and $\Sigma_t:=\big(\Sigma, F_t^* \overline g(t)\big)$, where $F(\,\cdot\,, t)$ is given by
\begin{enumerate}
\item[(i)] $F(\,\cdot\,,t):=\psi(\,\cdot\,,s(t))$, for $t \in (-\infty,T)$, in the shrinking case;
\smallskip
\item[(ii)] $F(\,\cdot\,,t):=\psi(\,\cdot\,,s(t))$, for $t \in (T,+\infty)$, in the expanding case; 
\smallskip
\item[(iii)] $F(\,\cdot\,,t):=\psi(\,\cdot\,,-t)$, for $t\in \mathbb{R}$, in the steady case.
\end{enumerate}

Assume that $\mathscr G:=\{\Sigma_t\}$ evolves by MCF in a gradient Ricci--Bourguignon soliton background on $M$. Then, it is a family of mean curvature solitons. 
Indeed, since $\overline g(t)=\sigma(t)\psi_t^*g$, we have $\nabla_g  f\circ \psi(\,\cdot\,,  s(t))= \sigma(t)\nabla_{\overline g(t)}  \overline f$, and then, for the shrinking and expanding cases, we obtain 
\begin{align*}
H(\,\cdot\,,t) &= \overline g(t)\Big( \frac{\partial}{\partial t} F(\,\cdot\,,t),  \nu(\,\cdot\,, t)\Big)= \overline g(t)\Big( \frac{\partial}{\partial s} \psi(\,\cdot\,,s(t)) \frac{ds}{dt},  \nu(\,\cdot\,, t)\Big)\\
 &=\overline g(t)\left(  \frac{ \nabla_g  f }{\sigma(s(t))}\circ \psi(\,\cdot\,,  s(t))\frac{\kappa \sigma'(t)}{\sigma^2(t)},  \nu(\,\cdot\,, t)\right)\\
 &=-\overline g(t)\left(  \frac{ \nabla_g  f}{\sigma(t)}\circ \psi(\,\cdot\,,  s(t)),  \nu(\,\cdot\,, t)\right)\\
 &=-\overline g(t)\Big(   \nabla_{\overline g(t)}  \overline f ,  \nu(\,\cdot\,, t)\Big)=- \nu(\,\cdot\,, t) \overline f.
\end{align*}
In the steady case, the proof is analogous, so the claim follows. 

It is worth noting that, as in case (iii), cases (i) and (ii) also yield $\psi_{s(t)}=\psi_t^{-1}$. Indeed, let $\{\Phi_\tau\}$ be the smooth one-parameter family of diffeomorphisms of $M$ generated by $\nabla_g f,$ with $\Phi_0=\operatorname{Id}.$ Since
\begin{equation*}
\frac{\partial}{\partial t}\psi_t
=\frac{1}{\sigma(t)}\nabla_g f\circ\psi_t,
\end{equation*}
the time reparametrization corresponding to each value of $\kappa$,
\begin{equation*}
a(t)=\int_{T-\kappa}^{t}\frac{du}{\sigma(u)}
=-\kappa\log\sigma(t)
\end{equation*}
gives
\begin{equation*}
\psi_t=\Phi_{-\kappa\log\sigma(t)} \quad\mbox{and}\quad \psi_t^{-1}=\Phi_{\kappa\log\sigma(t)}.
\end{equation*}
Now, note that $\sigma(s(t))=\sigma(t)^{-1}.$ Thus,
\begin{equation*}
\psi_{s(t)} =\Phi_{-\kappa\log\sigma(s(t))}= \Phi_{-\kappa\log\sigma(t)^{-1}}= \psi_t^{-1},
\end{equation*}
thereby proving our assertion.

\begin{proof}[\bf Proof of Theorem~\ref{NazasCharacterization}] 
Assume that $\Sigma$ is an $f$-minimal hypersurface of $(M,g)$, i.e., $H+\nu f=0$ on $\Sigma$. Let $\mathscr{G}=\{\Sigma_t\}$ be the smooth one-parameter family of isometric immersions of $\Sigma$ into $M$ defined as above. By conformal theory,  we obtain 
$H(\cdot, t)=\sigma(t)^{-1/2}H$ and $\nu(\cdot, t)\overline f = \sigma(t)^{-1/2}\nu f$, and then $H(\,\cdot\,,t)+\nu(\,\cdot\,,t)\,\overline{f}=0$. Hence,
\begin{align*}
\Big(\frac{\partial}{\partial t}  F (\,\cdot\,, t)\Big)^\perp &= \overline g(t)\Big(\frac{\partial}{\partial t}  F (\,\cdot\,, t), \nu(\,\cdot\,,t)\Big)\nu(\,\cdot\,,t)\\
&=- \nu(\,\cdot\,, t) (\overline f) \nu(\,\cdot\,, t)\\ 
&= H(\,\cdot\,,t) \nu(\,\cdot\,, t).   
\end{align*}
Now, we affirm that $F (\,\cdot\,, t)$ can be everywhere locally reparametrized to an MCF. Indeed, if 
\begin{equation*}
\frac{\partial}{\partial t} F(p,t)=H(p,t) \nu(p,t) + X(p,t),
\end{equation*}
where $X(p,t) \in d F_t(Tp\Sigma)\,\,\,\forall p \in \Sigma$. Let $\{\varphi_t\}$ be the smooth one-parameter family of diffeomorphisms of $\Sigma$ generated by $Y(p,t) =  -[dF_t]^{-1}(X(p,t))$, with $\varphi_0 = \operatorname{Id}_{\Sigma}$. Considering the reparametrization $\widetilde F(p,t) = F(\varphi_t (p), t)$, we have 
\begin{align*}
\frac{\partial}{\partial t}\widetilde F(p,t)=& dF(\varphi_{t}(p),t)\left(\frac{\partial}{\partial t}\varphi_{t}(p)\right)+\frac{\partial}{\partial t}F(\varphi_t(p),t) \\
=& dF(\varphi_{t}(p),t)\left(Y(\varphi_{t}(p),t)\right) + H(\varphi_t(p),t)\nu(\varphi_t(p),t)+X(\varphi_t(p),t)\\
=&H(\varphi_t(p),t)\nu(\varphi_t(p),t),
\end{align*}
which is enough to conclude that $\{\widetilde \Sigma _t := \widetilde F_t(\Sigma)\}$ evolves by MCF  in a gradient Ricci--Bourguignon soliton background on $M$. 

Now, we will show that any family $\mathscr F$ of mean curvature solitons is given by $\mathscr G$ up to reparametrization. Indeed, let $G_t$ be a smooth one-parameter family of isometric immersions of $\Sigma$ into $(M,\bar{g}(t))$, such that 
\begin{equation*}
\frac{\partial}{\partial t}G_t=\vec H_{G_t,\overline g(t)} \qquad\mbox{and}\qquad\vec H_{G_t,\overline g(t)} +\bigl(\nabla_{\overline g(t)}\overline f_t\bigr)^\perp=0,
\end{equation*}
where $\vec H_{G_t,\overline g(t)}$ stands for the mean curvature vector of $G_t$ in $(M,\overline g(t))$. Setting $\widehat G_t=\psi_t\circ G_t$, we have
\begin{equation*}
\frac{\partial}{\partial t} \widehat G_t =\frac1{\sigma(t)}\left(\nabla_gf\circ G_t + \vec H_{G_t,g}\right).
\end{equation*}
On the other hand
\begin{equation*}
\vec H_{G_t,g} +\bigl(\nabla_gf\circ G_t\bigr)^\perp=0,
\end{equation*}
and therefore
\begin{equation*}
\frac{\partial}{\partial t}\widehat G_t =\frac1{\sigma(t)}\bigl(\nabla_gf\circ\widehat G_t\bigr)^\top.
\end{equation*}
However, as we have seen above, we can reparametrize it to obtain a map $\widehat G_t$ that is independent of $t$. Consequently, $G_t$ is, up to reparametrization, the family obtained by applying $\psi_t^{-1}$ to a fixed $f$-minimal immersion.
\end{proof} 

The proofs of Theorems~\ref{Cor-FTS} and~\ref{Huisken_monotonicity} rely on the following two lemmas.

\begin{lemma}\label{evolução do RB}
For any MCF in a Ricci--Bourguignon flow background, the following equation holds on $\Sigma_t$
\begin{align}
\frac{\partial}{\partial t}g_{ij} = -2R_{ij} + 2\rho R g_{ij} - 2H\mathcal{A}_{ij}.
\end{align}
\end{lemma}

\begin{proof}
On $\Sigma_t$ we obtain
\begin{align*}
\frac{\partial}{\partial t}g_{ij} :=& \frac{\partial}{\partial t}F_{t}^{*}g_{ij}\\
=& \frac{\partial}{\partial t}g\left(\frac{\partial F_{t}}{\partial x_{i}}, \frac{\partial F_{t}}{\partial x_{j}}\right) + g\left(\nabla_{\frac{\partial F_{t}}{\partial t}}\frac{\partial F_{t}}{\partial x_{i}}, \frac{\partial F_{t}}{\partial x_{j}}\right) + g\left(\frac{\partial F_{t}}{\partial x_{i}},\nabla_{\frac{\partial F_{t}}{\partial t}}\frac{\partial F_{t}}{\partial x_{j}}\right)\\
=& -2R_{ij} + 2\rho R g_{ij} + g\left(\nabla_{\frac{\partial F_{t}}{\partial t}}\frac{\partial F_{t}}{\partial x_{i}}, \frac{\partial F_{t}}{\partial x_{j}}\right) + g\left(\nabla_{\frac{\partial F_{t}}{\partial t}}\frac{\partial F_{t}}{\partial x_{j}}, \frac{\partial F_{t}}{\partial x_{i}}\right).
\end{align*}
On the other hand, we have
\begin{align*}
g\left(\nabla_{\frac{\partial F_{t}}{\partial t}}\frac{\partial F_{t}}{\partial x_{i}}, \frac{\partial F_{t}}{\partial x_{j}}\right)\! = \! Hg\left(\nabla_{\nu_t}\frac{\partial F_{t}}{\partial x_{i}}, \frac{\partial F_{t}}{\partial x_{j}}\right) \!=\! -Hg\left(\nabla_{\frac{\partial F_{t}}{\partial x_{j}}} \frac{\partial F_{t}}{\partial x_{i}}, \nu_t\right)\! :=\! -H\mathcal{A}_{ji}
\end{align*}
which is enough to conclude the result of the lemma.
\end{proof}
\begin{lemma}\label{Huisken:lem}
For any MCF in a Ricci--Bourguignon flow background, the following equation holds on $\Sigma_t$
\begin{align*}
\dfrac{d}{dt}(d A_{g}) = -\left(R^{i}\!_{i} - \rho(n-1)R + H^{2}_{g}\right)dA_{g}.
\end{align*}
\end{lemma}
\begin{proof}
From the well-known area formula and Lemma~\eqref{evolução do RB}, we have
\begin{align*}
\dfrac{d}{dt}(d A_{g}) =& \dfrac{1}{2}g^{ij}\dfrac{\partial}{\partial t} g_{ij} dA_{g}\\
=& \dfrac{1}{2}g^{ij}\left(  -2R_{ij} + 2\rho R g_{ij} - 2H_{g} \mathcal{A}_{ij}\right)dA_{g}\\
=& -\left(R^{i}\!_{i} - \rho(n-1)R + H^{2}_{g}\right)dA_{g}
\end{align*}
thereby proving the lemma.
\end{proof}
\begin{proof}[\textbf{Proof of Theorem~\ref{Cor-FTS}}] 
Taking $n=3$ in Lemma~\ref{Huisken:lem}, we have
\begin{align}\label{area}
\frac{d}{dt}\operatorname{Area}(\Sigma_{t}) = -\int_{\Sigma_t}\left(R^{i}\!_{i} - 2\rho R + H_{g}^{2}\right)dA_{g}. 
\end{align} 
Choose a point $p \in \Sigma_t$ and let $\{e_1, e_2, e_3\}$ be an orthonormal basis for $T_pM$ such that $e_1$ and $e_2$ are tangent to $\Sigma_t$. Denoting the sectional curvatures as 
\begin{equation*} 
K(e_{i},e_{j})=R(e_{i},e_{j},e_{i},e_{j}), \quad i,j \in \{1,2,3\},
\end{equation*} 
we have $R=2[K(e_{1},e_{2})+K(e_{1},e_{3})+K(e_{2},e_{3})]$ and $R_{33}=K(e_{1},e_{3})+K(e_{2},e_{3})$. Then, 
\begin{align*} 
R^{i}\!_{i} &= \sum_{i=1}^{2}R_{ii} = R - R_{33} = R - (K(e_{1},e_{3})+K(e_{2},e_{3}))\\
&= R - \left(\frac{R}{2}-K(e_{1},e_{2})\right) = \frac{1}{2}(R+R^{ij}\!_{ij}),    
\end{align*} 
where in the last equality we have used
\begin{equation*} 
R^{ij}\!_{ij} = g^{im}g^{jk}R_{mkij}=R_{1212}+R_{2121}=2K(e_{1},e_{2}). 
\end{equation*} 
By Gauss's formula, we obtain
\begin{equation}\label{traçoparcial}
R^{i}\!_{i}=\sum_{i=1}^{2}R_{ii} = \frac{1}{2}\left(R + R^{ij}\!_{ij} \right) = \frac{1}{2}\left(R + \widehat{R} - H_{g}^2 + \mathcal{A}^{ij}\mathcal{A}_{ij} \right), 
\end{equation} 
where $\widehat{R}$ denotes the scalar curvature of $\Sigma_t$. Replacing \eqref{traçoparcial} into \eqref{area}, one has
\begin{equation}\label{derivada}
\frac{d}{dt}\operatorname{Area}(\Sigma_{t}) = -\frac{1}{2}\int_{\Sigma_t}\left((1-4\rho)R + \widehat{R} + H_{g}^{2} + \mathcal{A}^{ij}\mathcal{A}_{ij} \right)dA_{g}.
\end{equation}
Since $M$ is compact, Proposition~4.1 in~\cite{Catino} ensures that if $R(\cdot,0)\geq\alpha$ for some $\alpha \in \mathbb{R}$ then
\begin{equation*}
R_{\min}(t)\geq\frac{3\alpha}{3-2(1-3\rho)\alpha t}.
\end{equation*}
Note that we can choose $\beta>0$ such that $R(\cdot,0)\geq-\beta$. So, from the previous estimate, one has
\begin{equation}\label{estima0}
 R(\cdot,t)\geq
-\frac{3\beta}{3+2(1-3\rho)\beta t}.
\end{equation}
Since $1-3\rho>0$, it follows from \eqref{estima0} that, for \(t>0\),
\begin{equation}\label{estima}
 R(\cdot,t)\geq-\frac{3}{2(1-3\rho)t}.
\end{equation}
Moreover, by assumption $1-4\rho\geq0$, then \eqref{derivada}, \eqref{estima} and Gauss--Bonnet theorem imply
\begin{equation}\label{estima1}
\frac{d}{dt}\operatorname{Area}(\Sigma_{t}) \leq\frac{3(1-4\rho)}{4(1-3\rho)}
\frac{\operatorname{Area}(\Sigma_{t})}t-4\pi.
\end{equation}
Defining $\tau =\frac{3(1-4\rho)}{4(1-3\rho)}$, we have $0\leq\tau<1$ and $1-\tau=\frac1{4(1-3\rho)}>0$.
Multiplying \eqref{estima1} by the integrating factor $t^{-\tau}$, we obtain
\begin{equation*}
\frac{d}{dt}\left(t^{-\tau}\operatorname{Area}(\Sigma_t)\right)\leq -4\pi t^{-\tau}.
\end{equation*}
Fixing any \(t_0\in(0,S)\) and integrating from \(t_0\) to \(t\), we find
\begin{equation*}
t^{-\tau}\operatorname{Area}(\Sigma_t)\leq t_0^{-\tau}\operatorname{Area}(\Sigma_{t_0})-\frac{4\pi}{1-\tau}\left(t^{1-\tau}-t_0^{1-\tau}\right).
\end{equation*}
Since
\begin{equation*}
1-\tau=\frac{1}{4(1-3\rho)}>0,
\end{equation*}
it follows that, for any \(t_0>0\) such that the MCF exists at time \(t_0\), we would have \(\operatorname{Area}(\Sigma_t)\leq 0\) for
\begin{equation*}
t\geq t_0\left(1+\frac{\operatorname{Area}(\Sigma_{t_0})}{16\pi(1-3\rho)t_0}\right)^{4(1-3\rho)}.
\end{equation*}
Therefore, the maximal interval of smooth existence of the MCF within the given ambient flow is finite.
\end{proof}

We finish this section by showing Theorem~\ref{Huisken_monotonicity}.

\begin{proof}[\bf Proof of Theorem~\ref{Huisken_monotonicity}]
Lemma~\ref{Huisken:lem} and a straightforward computation  yield
\begin{align*}
\dfrac{d}{dt}\int_{\Sigma_t} e^{-\overline{f}} d A_{\overline{g}} &= - \int_{\Sigma_t} \Big(\frac{d}{dt}\overline{f} + \overline{R}^i\!_i + H_{\overline{g}}^2 -\rho(n-1)\overline{R} \Big) e^{-\overline{f} } d A_{\overline{g}}.
\end{align*}
Therefore,
\begin{align*}
\dfrac{d}{dt}\int_{\Sigma_t} e^{-\overline{f}} d A_{\overline{g}}  &= - \int_{\Sigma_t} \Big(\dfrac{\partial}{\partial t}\overline{f} + H_{\overline{g}} \nu_t \overline{f}  + \overline{R}^i\!_i + H_{\overline{g}}^2  -\rho(n-1)\overline{R}\Big) e^{-\overline{f}} d A_{\overline{g}}.
\end{align*}

First, assume $\overline{g}(t)$ is a gradient steady soliton. In this case, we can take the trace of equation ~\eqref{mod_grad_Ricci_soliton}  on $\Sigma_t$ to get
\begin{align*}
0=\overline{R}^i\!_i + \overline{\nabla}^i \overline{\nabla}_i \overline{f} -\rho(n-1)\overline{R} = \overline{R}^i\!_i + \widehat{\nabla}^i \widehat{\nabla}_i \overline{f} - H_{\overline{g}} \nu_t \overline{f}  -\rho(n-1)\overline{R}.
\end{align*}
Then, using \eqref{self-solution}, we obtain
\begin{align*}
\dfrac{d}{dt}\int_{\Sigma_t} e^{-\overline{f}} d A_{\overline{g}}
&= - \int_{\Sigma_t} \Big(|\nabla_{\overline{g}} \overline{f}|_{\overline{g}}^2 - \widehat{\Delta}_{\overline{g}}\overline{f} +2 H_{\overline{g}} \nu_t \overline{f}  + H_{\overline{g}}^2 \Big)e^{-\overline{f}} d A_{\overline{g}} \\
&= - \int_{\Sigma_t} \Big( |\widehat{\nabla}_{\overline{g}} \overline{f}|_{\overline{g}}^2 +(\nu_t \overline{f})^2 - \widehat{\Delta}_{\overline{g}}\overline{f} + 2H_{\overline{g}} \nu_t \overline{f}  + H_{\overline{g}}^2\Big)e^{-\overline{f}} d A_{\overline{g}}\\
&= - \int_{\Sigma_t} \Big(H_{\overline{g}} +\nu_t \overline{f} \Big)^2 e^{-\overline{f}} d A_{\overline{g}},
\end{align*}
where in the second line we have used the equality
\begin{equation*}
\widehat{\Delta}_{\overline{g}} e^{- \overline{f}} = (|\widehat{\nabla}_{\overline{g}} \overline{f}|_{\overline{g}}^2 -  \widehat{\Delta}_{\overline{g}}\overline{f})e^{-\overline{f}}
\end{equation*}
and Stokes' theorem. Since the boundary integrand on the right-hand side is nonnegative, we have the result of the theorem for the steady case.

For the shrinking case, we claim that the function
\begin{equation*}
(-\infty,T)\ni t\mapsto [4\pi(T-t)]^{\frac{-(n-1)}{2}}\int_{\Sigma_t} e^{-\overline{f}} dA_{\overline{g}}    
\end{equation*}
is non-increasing. Indeed, as above, we take the
trace of equation~\eqref{mod_grad_Ricci_soliton} on $\Sigma_t$ to obtain
\begin{align*}
\dfrac{n - 1}{2(T-t)}=\overline{R}^i\!_i + 
\overline{\nabla}^i \overline{\nabla}_i \overline{f}  -\rho(n-1)\overline{R} 
= \overline{R}^i\!_i + \widehat{\nabla}^i \widehat{\nabla}_i \overline{f}- 
H_{\overline{g}} \nu_t \overline{f}  -\rho(n-1)\overline{R}.
\end{align*}
Then,
\begin{align*}
&\dfrac{d}{dt}\Big([4\pi(T-t)]^{\frac{-(n-1)}{2}}\int_{\Sigma_t} e^{-\overline{f}} d A_{\overline{g}}\Big) = \frac{d}{dt}\left([4\pi(T-t)]^{\frac{-(n-1)}{2}}\right)\int_{\Sigma_t}e^{-\overline{f}} dA_{\overline{g}}\\
&- [4\pi(T-t)]^{\frac{-(n-1)}{2}}\int_{\Sigma_t}\dfrac{n - 1}{2(T-t)}e^{-\overline{f}}dA_{\overline{g}}\\
&- [4\pi(T-t)]^{\frac{-(n-1)}{2}}\int_{\Sigma_t}\Big( |\widehat{\nabla}_{\overline{g}} \overline{f}|_{\overline{g}}^2 +(\nu_t \overline{f})^2 - \widehat{\Delta}_{\overline{g}}\overline{f} + 2H_{\overline{g}} \nu_t \overline{f}  + H_{\overline{g}}^2\Big)e^{-\overline{f}} d A_{\overline{g}}\\
&= -[4\pi(T-t)]^{\frac{-(n-1)}{2}} \int_{\Sigma_t} \Big(H_{\overline{g}} +\nu_t \overline{f} \Big)^2 e^{-\overline{f}} dA_{\overline{g}}.
\end{align*}
This proves the claim, and so the theorem for the shrinking case. Finally, in a similar way, one proves the expanding case.
\end{proof}

\section{\bf An application of the Huisken monotonicity-type formula}\label{sect:apllication}
In this section, we prove Theorem~\ref{application:Huisken}, which is an application of a Huisken monotonicity-type formula for the context of MCF in a gradient Ricci--Bourguignon soliton background. We begin by providing a more detailed description of the normalized family $\widetilde{\mathscr{F}}$ of the MCF in a gradient shrinking Ricci--Bourguignon soliton background.

Let $\mathscr F :=\{\Sigma_t; \; t\in [0,T)\}$ be the MCF in a gradient shrinking Ricci--Bourguignon soliton background. Write $H_{\overline g}=H(p,t),$  $\nu_t=\nu(p,t)$ and $\widetilde F_{s} = \psi_t\circ F_t,$ for short. Setting $s=-\log (T- t)$ and $\widetilde{\Sigma}_s = (\Sigma,\widetilde F_{s}^{*}g),$ one has $s\in [-\log T,\infty)$, $\frac{ds}{dt}=\frac{1}{T- t}$ and $\frac{dt}{ds}=T- t.$ Since $\overline g(t)= (T- t)\psi^*_t g$, we have $H_{\overline{g}} =\frac{1}{\sqrt{T- t}} H_{\psi^*_t g}$ and $\nu_t = \frac{1}{\sqrt{T- t}}\nu_{\psi^*_t g}$. Thus, we recall the family $\widetilde{\mathscr F}:=\{\widetilde \Sigma_s;\; s\in [-\log T,\infty)\}$ which is the normalized MCF of $\mathscr F$ and satisfies
\begin{align*}
\frac{\partial \widetilde F_s}{\partial s}=&  \left[ \frac{d \psi_t}{dt} \circ F_t +(\psi_t)_* \left(\frac{\partial F_t}{\partial t}\right)\right]\frac{dt}{ds}\\
=&  \frac{\nabla_{g}f\circ\psi_t}{T - t}  \circ F_t \cdot \frac{dt}{ds} + (T-  t)(\psi_t)_* \big(H_{\overline{g}}\nu_t\big) \\
=& \big(\nabla_g f + H_g\big)(\widetilde F_s) .
\end{align*}

We now turn to the proof of our next main result, which is by means of Theorem~\ref{Huisken_monotonicity} and of the Arzelà--Ascoli theorem. The latter provides a sequence of isometric immersions on an exhaustion of $\Sigma_\infty$ converging to a limiting global solution which is an $f_\infty$-minimal hypersurface of $(M, g)$. We prove global supremum estimates depending only on the initial bounds on the Riemann tensor $\operatorname{Rm}_g$, whereas the interior estimates depend on the full $C^\infty$ norm of $\overline g$ (see Section~\ref{sect:apllication}). In what follows, we briefly lay the groundwork for the compact and noncompact cases.

\begin{definition}[Bounded geometry]\label{DefBGeom}
We say that a complete Riemannian manifold $(M,g)$ has bounded geometry if for every integer $j\geq 0$ there exist positive constants $C_j$ and $\eta$ such that the Riemann curvature tensor $\operatorname{Rm}_g$ and the injectivity radius $\operatorname{inj}(M, g)$ satisfy
\begin{equation*}
|\nabla_g^j\operatorname{Rm}_g| < C_j  \quad \mbox{and} \quad \operatorname{inj}(M, g) \geq \eta>0.
\end{equation*}
\end{definition}

Notice that $\mathbb{R}^n$ with the standard metric and compact Riemannian manifolds have bounded geometry. Actually, bounded geometry provides a natural generalization to the known settings of compact and Euclidean spaces. In particular, the class of manifolds of bounded geometry allows us to uniformly apply constructions that are well known for compact smooth manifolds. For a way to construct complete manifolds of bounded geometry, we refer to Eldering~\cite{eldering2013normally}.

The notion of the convergence of immersions from pointed manifolds is defined as follows. It is the immersion map version of the Cheeger--Gromov convergence.

\begin{definition}\label{pointed:manifold}
Let $(M, g)$ be an  $n$-dimensional complete Riemannian manifold with bounded geometry. Assume that for each $k \geq 1$ we have an $m$-dimensional pointed manifold $\left(\Sigma_{k}, p_{k}\right)$ and an immersion $F_{k}: \Sigma_{k} \rightarrow M$. Then, we say that a sequence of immersions $\left\{F_{k}: \Sigma_{k} \rightarrow M\right\}_{k=1}^{\infty}$ converges to an immersion $F_{\infty}: \Sigma_{\infty} \rightarrow M$ from an $m$-dimensional pointed manifold $\left(\Sigma_{\infty}, p_{\infty}\right)$ if there exist
\begin{enumerate}
\item An exhaustion $\left\{U_{k}\right\}_{k=1}^{\infty}$ of $\Sigma_{\infty}$ with $p_{\infty} \in U_{k}$, and
\smallskip
\item A sequence of diffeomorphisms $\Psi_{k}: U_{k} \rightarrow V_{k} \subset \Sigma_{k}$ with $\Psi_{k}\left(p_{\infty}\right)=p_{k}$ such that the sequence of maps $F_{k} \circ \Psi_{k}: U_{k} \rightarrow M$ converges in $C^{\infty}$ to $F_{\infty}: \Sigma_{\infty} \rightarrow M$ uniformly on compact sets in $\Sigma_{\infty}$. 
\end{enumerate} 
\end{definition} 

To apply the Arzelà--Ascoli theorem, we need a uniform interior estimate for $\big|\widehat\nabla_g^k \mathcal A(\widetilde F_s)\big|$ on $\widetilde\Sigma_s\in\widetilde{\mathscr F}$, analogous to those established by Ecker and Huisken~\cite[Sect.3]{ecker1991interior} for MCF in Euclidean space, by Yamamoto~\cite[Sect.4]{yamamoto2020meancurvature} for MCF in a Ricci flow background, and by the first and second authors and Yamamoto~\cite{gomes2024anote} for MCF in an extended Ricci soliton background.

\begin{proposition}\label{Underthesamesatupofproposition01}
Let $(M,g)$ be an $(n\geq 3)$-dimensional complete Riemannian manifold with bounded geometry, and let $\mathscr F:=\{\Sigma_t; \; t\in [0,T)\}$ be the MCF of $\Sigma$ in a gradient shrinking Ricci--Bourguignon soliton background which develops a singularity of type-I. Consider $\widetilde{\mathscr F}:=\{\widetilde \Sigma_s;\; s\in [-\log T,\infty)\}$ in $(M,g)$. Then, for every integer $k\geq0$ there is a positive constant $C_k$ which does not depend on $s$ such that
\begin{equation*}
\big|\widehat\nabla_g^k \mathcal A(\widetilde F_s)\big|_g \leq C_k \quad \mbox{on} \quad  \Sigma \times [-\log T, \infty),  
\end{equation*}
where $\widehat\nabla_g$ is defined from the Levi–Civita connection on $\big(\Sigma, \widetilde F^*_sg\big).$
\end{proposition}
\begin{proof}
The proof follows the argument of Yamamoto, since metric compatibility implies that every covariant derivative of the scalar curvature is a contraction of the corresponding covariant derivative of the ambient curvature tensor. Consequently, the additional Ricci--Bourguignon terms belong to the same classes of $*$-products and have the same scaling degrees as those appearing in Yamamoto’s estimates, so the induction applies with constants additionally depending on $\rho$ (see~\cite[Prop.~4.9]{yamamoto2020meancurvature}).
\end{proof}

The following lemma gives the variation of the weighted area functional on the normalized MCF.
\begin{lemma}\label{Prop:4.6}
Let $(M,g)$ be an $(n\geq 3)$-dimensional Riemannian manifold, and let $\mathscr F:=\{\Sigma_t; \; t\in [0,T)\}$ be the MCF of $\Sigma$ in a gradient shrinking Ricci--Bourguignon soliton background. Consider  $\widetilde{\mathscr F}:=\{\widetilde \Sigma_s;\; s\in [-\log T,\infty)\}$ in $(M,g)$. Then, for $s \in [-\log T , \infty)$, we have
\begin{equation*}
\frac{d}{ds}\int_\Sigma e^{-f \circ \widetilde F_s}dA_{\widetilde F_s^*g} = -\int_\Sigma \big( H_g(\widetilde F_s) + \nu_s  (f \circ \widetilde F_s) \big)^2  e^{-f \circ \widetilde F_s}dA_{\widetilde F_s^*g}.
\end{equation*}
\end{lemma}

\begin{proof}
As $\overline f\circ F_t =  f \circ \widetilde F_s$ and $F^*_t \overline g= (T-t)  (\psi_t \circ F_t)^* g = (T-t)\widetilde F^*_s g$ both on $\Sigma$, we have
\begin{align}
[4\pi(T - t)]^{-\frac{n - 1}{2}} e^{-\overline f} dA_{\overline g}&=(4\pi)^{-\frac{n -1}{2}} e^{-f\circ \widetilde F_s} dA_{\widetilde F_s^*g},\label{xa:era:01}\\
(T - t)\left(H_{\overline g}+ \nu_t \overline  f \right)^2 &= \left(H_g(\widetilde F_s) + \nu_s (f \circ \widetilde F_s) \right)^2.\label{xa:era:02} 
\end{align}
The result of the lemma follows from Theorem~\ref{Huisken_monotonicity}, identities \eqref{xa:era:01} and \eqref{xa:era:02}, together with the chain rule.
\end{proof}

The following lemmas, which rely on a Hamilton-type identity, are also needed to prove Theorem~\ref{application:Huisken}.

First, assume that $\rho\neq 0$ in equation~\eqref{Bizu1} of Section~\ref{example:section}. Following the approach in Hamilton's identity (see~Hamilton~\cite{Hamilton} or, alternatively, Cao~\cite{HDCao}), we obtain
\begin{equation*}
(1-2(n-1)\rho)R_{g}+|\nabla_{g} f|^{2}-2\lambda f - 2\rho \phi = C, 
\end{equation*}
where $\phi$ is given by $d\phi=Rdf$, and $C$ is a constant. In the shrinking case, we can assume, after rescaling the metric $g$ and adding a constant to $f$ if necessary, that the previous equation takes the form 
\begin{align*}
(1-2(n-1)\rho)R_{g} + |\nabla_{g} f|^2 - f -2\rho\phi = 0.
\end{align*}
Since $\bar g(t)=\sigma(t)\psi_t^*g,$  where $\sigma(t)=T-t$, then 
\begin{align*}
(1-2(n-1)\rho)R_{\psi_{t}^{*}g} + |\nabla_{\psi_{t}^{*}g} \overline{f}|^2 - \overline{f} - 2\rho \overline{\phi}= 0,
\end{align*}
where $\overline{f} = \psi^{*}_{t}f$ and $\overline{\phi} = \psi^{*}_{t}\phi$. By conformal theory, we get $R_{\bar{g}}=R_{\psi_t^*g}/\sigma(t)$ and $|\nabla_{\bar g} \bar f|^2=|\nabla_{\psi_t^*g} \bar f|^2/\sigma(t)$, whence we obtain the Hamilton-type identity
\begin{equation}\label{like Hamilton Identity}
(1-2(n-1)\rho)R_{\overline{g}} + |\nabla_{\overline{g}}\overline{f}|^{2}_{\overline{g}}-\frac{\overline{f}+2\rho\overline{\phi}}{T-t}=0 
\end{equation}
along the gradient Ricci--Bourguignon soliton on $M$.

We also note the following facts for later use. The scalar curvature of a complete ancient solution to the Ricci--Bourguignon flow on a compact manifold is nonnegative along the flow (see Catino et al.~\cite[Prop.~4.4]{Catino}). Likewise, the scalar curvature of a complete ancient solution to the Schouten flow on a noncompact manifold is nonnegative along the flow (see Catino and Mazzieri~\cite[Prop.~5.1]{catino2016gradient}). Since gradient shrinking soliton solutions are ancient, this nonnegativity condition is automatically satisfied in the compact setting of Theorem~\ref{application:Huisken} and in the noncompact setting of Theorem~\ref{Thmnoncommpact}. Consequently, the assumption that the scalar curvature remains nonnegative along the flow, imposed in Lemmas~\ref{prop:4.7} and \ref{Lem:4.10} below, does not yield additional restriction in either of the settings in which these lemmas will be applied.

\begin{lemma}\label{prop:4.7}
Let $(M,g)$ be an $(n\geq 3)$-dimensional complete Riemannian manifold with bounded geometry, and let $\mathscr F:=\{\Sigma_t; \; t\in [0,T)\}$ be the MCF of $\Sigma$ in a gradient shrinking Ricci--Bourguignon soliton background which develops a singularity of type-I. Consider $\widetilde{\mathscr F}:=\{\widetilde \Sigma_s;\; s\in [-\log T,\infty)\}$ in $(M,g)$ so that $\phi\circ \widetilde{F}_{s},$ with $\rho \leq 1/2(n-1)$ and $\rho\neq0,$ remains uniformly bounded on $[-\log T,\infty).$ In addition, assume that the scalar curvature $R_{\overline g(t)}$ is nonnegative along the flow. Then, there exists a positive constant $C$ such that
\begin{equation*}
\int_\Sigma e^{-\frac{f}{2} \circ \widetilde F_s} dA_{\widetilde F^*_s g}\leq C   
\end{equation*}
uniformly on $[-\log T, \infty).$
\end{lemma}
\begin{proof}
We begin by substituting $\overline{\alpha}_{t}:=\frac{1}{[4\pi(T-t)]^{\frac{n}{2}}}e^{-\frac{\overline f}{2}}$, $u_{t}:=[4\pi(T-t)]^{1/2},$ $h:=-2\operatorname{Ric}_{\overline g}+2\rho R_{\overline{g}}\overline{g}$ and $V:= H_{\overline g}\nu_t$ into Prop.~3.2 of \cite{yamamoto2020meancurvature} to obtain 
\begin{eqnarray*}
&&\frac{d}{dt} \int_{\Sigma}u_{t}\mathop{F_{t}^{*}\hspace{-0.5mm}\overline{\alpha}_{t}} dA_{\overline g}\\
&=&-\int_{\Sigma}u_{t}\Bigl(H_{\overline g}+\frac{1}{2} \nu_t\overline f\Bigr)^2\mathop{F_{t}^{*}\hspace{-0.5mm}\overline{\alpha}_{t}} dA_{\overline g}\\
&&+\int_{\Sigma}u_{t} F_{t}^{*}\left(\Delta_{\overline g}\overline \alpha_{t}+\frac{\partial \overline \alpha_t}{\partial t}-(1-n\rho)\overline\alpha_tR_{\overline g}\right)  dA_{\overline g}\\
&& +\int_{\Sigma}\left( \frac{\partial u_{t}}{\partial t} -\widehat\Delta_{\overline g}u_{t}+
u_{t}\left(\frac{1}{2}\mathrm{Hess}_{\overline g}\overline f-\frac{h}{2}\right)(\nu_{t},\nu_{t})\right)\mathop{F_{t}^{*}\hspace{-0.5mm}\overline{\alpha}_{t}} dA_{\overline g}. 
\end{eqnarray*}
By using $\Delta_{\overline g} \overline f=-(1-n\rho)R_{\overline g}+\frac{n}{2(T-t)}$ (see \eqref{mod_grad_Ricci_soliton}), $R_{\overline g} \geq 0$ \eqref{self-solution} and $|\nabla \overline f|^{2} =\frac{\overline{f}+2\rho\overline{\phi}}{T-t} - (1-2(n-1)\rho)R_{\overline{g}}$ (see \eqref{like Hamilton Identity}), we obtain 

\begin{eqnarray*}
\Delta_{\overline g}\overline \alpha_{t}+\frac{\partial \overline \alpha_t}{\partial t}-(1-n\rho)\overline\alpha_tR_{\overline g} &=&{\overline\alpha_{t}}\biggl(\!-\dfrac{\overline{f}+2\rho\overline{\phi}}{4(T-t)}-\dfrac{(1-2\rho)R_{\overline g}}{4}+\dfrac{n}{4(T-t)}\!\! \biggr)\\
 &&\leq\dfrac{\overline \alpha_t}{4(T-t)}\left(n-(\overline{f}+2\rho\overline{\phi})\right).
\end{eqnarray*}
Furthermore, a simple calculation shows that $u_t$ satisfies 
\begin{equation*}
\frac{\partial u_{t}}{\partial t} -\widehat\Delta_{\overline g}u_{t}+u_{t}\left(\mathrm{Hess}_{\overline g}\overline f-\frac{h}{2}\right)(\nu_{t},\nu_{t})=0,
\end{equation*}
which is equivalent to
\begin{equation*}
\frac{\partial u_{t}}{\partial t} -\widehat\Delta_{\overline g}u_{t}+u_{t}\left(\frac{1}{2}\mathrm{Hess}_{\overline g}\overline f-\frac{h}{2}\right)(\nu_{t},\nu_{t})=-\frac{1}{2}u_{t}\mathrm{Hess}_{\overline g}\overline f(\nu_{t},\nu_{t}).
\end{equation*}
Equation $\mathrm{Hess}_{\overline g}\overline f=\frac{1}{2(T-t)}\overline g-\mathrm{Ric}_{\overline g} + \rho R_{\overline g}\overline g$ implies 
\begin{equation*}
-\frac{1}{2}u_{t}\mathrm{Hess}_{\overline g}\overline f(\nu_{t},\nu_{t})=u_{t}\left(-\frac{1}{4(T-t)}+\frac{1}{2}(\mathrm{Ric}_{\overline g} - \rho R_{\overline{g}}\overline{g})(\nu_{t},\nu_{t})\right).
\end{equation*}
On the other hand, since $(M,g)$ has bounded geometry, we have that $C'':=\max_{M}\big\{|\mathrm{Ric}_g|_g + |\rho R_{g}g|_g\big\}$ is a constant, then 
\begin{eqnarray*}
&&(\mathrm{Ric}_{\overline g}-\rho R_{\overline{g}}\overline{g})(\nu_{t},\nu_{t})\\
&\leq& |\mathrm{Ric}_{\overline g} -\rho R_{\overline{g}}\overline{g}|_{\overline g}\leq\frac{|\mathrm{Ric}_g|_g + |\rho R_g g|_g}{T-t}\leq \frac{C''}{T-t}.    
\end{eqnarray*}
Hence,
\begin{equation*}
\frac{d}{dt} \int_{\Sigma}u_{t}\mathop{F_{t}^{*}\hspace{-0.5mm}\overline{\alpha}_{t}} dA_{\overline g}
\leq \frac{1}{4(T-t)}\int_{\Sigma}\Bigl(C_{0}-(\overline{f}+2\rho\overline{\phi})\circ F_{t}\Bigr)u_{t} \mathop{F_{t}^{*}\hspace{-0.5mm}\overline{\alpha}_{t}} dA_{\overline g}, 
\end{equation*}
where $C_{0}:=n+2C''$. Now, set $W(s)=\int_\Sigma e^{-\frac{f}{2}\circ \widetilde F_{s}}dA_{\widetilde F_s^*g},$ for $s=-\log(T-t)$. Thus, combining a straightforward computation with conformal theory and the previous estimate, we obtain
\begin{align*}
W'(s)    =&(4\pi)^{\frac{n-1}{2}}(T-t)\frac{d}{dt}\int_\Sigma u_{t}\mathop{F_{t}^{*}\hspace{-0.5mm}\overline{\alpha}_{t}} dA_{\overline g}\\
\leq&\frac14\int_\Sigma
\bigl(C_0-(f+2\rho\phi)\circ \widetilde F_{s}\bigr)e^{-\frac{f}{2}\circ \widetilde F_{s}}dA_{\widetilde F_s^*g}.
\end{align*}
Next, by Lemma~\ref{Prop:4.6} and the assumption on $\phi\circ\widetilde F_s,$  we have
\begin{equation*}
J(s):=\int_\Sigma e^{-f\circ \widetilde F_{s}}dA_{\widetilde F_s^*g}\leq J(s_0),
\quad \mbox{and}\quad B:=\sup_{\Sigma\times[-\log T, \infty) }|2\rho\phi\circ\widetilde F_s|,
\end{equation*}
for some constants $J(s_0)$ and $B$. For $D:=C_0+B$ and every real $x$, it is true that
\begin{equation*}
(D-x)e^{-x/2}
\leq-e^{-x/2}+2e^{(D-1)/2}e^{-x},
\end{equation*}
which follows from the fact that the function $q(x)=(D+1-x)e^{x/2}$ attains its global maximum at $x=D-1$, where it takes the value $2e^{(D-1)/2}$. Thus,
\begin{equation*}
W'(s)\leq-\frac14W(s)
+\frac12e^{(D-1)/2}J(s_0).
\end{equation*}
Multiplying by $e^{(s-s_0)/4}$ and integrating over $[s_0,s]$, we obtain
\begin{equation*}
W(s)\leq e^{-(s-s_0)/4}W(s_0)+2e^{(D-1)/2}J(s_0)\bigl(1-e^{-(s-s_0)/4}\bigr).
\end{equation*}
In particular,
\begin{equation*}
\sup_{s\geq s_0}W(s) \leq \max\left\{W(s_0),\,2e^{(D-1)/2}J(s_0)\right\}<\infty.
\end{equation*}
This yields the uniform bound asserted in the lemma.
\end{proof}

\begin{lemma}\label{Lem:4.10}
Let $(M,g)$ be an $(n\geq 3)$-dimensional complete Riemannian manifold with bounded geometry, and let $\mathscr F:=\{\Sigma_t; \; t\in [0, T)\}$ be the MCF of $\Sigma$ in a gradient shrinking Ricci--Bourguignon soliton background which develops a singularity of type-I. Consider $\widetilde{\mathscr F}:=\{\widetilde \Sigma_s;\; s\in [-\log T,\infty)\}$ in $(M,g)$ such that $\phi\circ \widetilde{F}_{s},$ with $\rho \leq 1/2(n-1)$ and $\rho\neq0,$ remains uniformly bounded on $[-\log T,\infty).$ In addition, assume that the scalar curvature $R_{\overline g(t)}$ is nonnegative along the flow. Then, there exists a positive constant $C'$ such that
\begin{equation*}
\left|\frac{d^2}{d^2s}\int_\Sigma e^{-f\circ \widetilde F_s}dA_{\widetilde F^*_s g}\right|  = \left|\frac{d}{ds}\int_\Sigma \big(H_g(\widetilde F_s) + \nu_s  (f \circ \widetilde F_s) \big)^2e^{-f\circ \widetilde F_s}dA_{\widetilde F_s^*g}\right|\leq C' 
\end{equation*}
uniformly on $[-\log T, \infty).$ 
\end{lemma}
\begin{proof}
We already know that we can assume $(1-2(n-1)\rho)R+|\nabla f|^2-f-2\rho\phi=0$ and $R\geq0$ along the gradient shrinking Ricci--Bourguignon soliton on $M$. So, 
$0\leq |\nabla f|^2 \leq f+2\rho \phi$. The lemma follows from Lemma~\ref{prop:4.7} and by following the same arguments as in Yamamoto~\cite{yamamoto2020meancurvature}.
\end{proof}
Before proceeding to the proof of one of the main results of this paper. We show that, given a sequence of isometric immersions in $\widetilde{\mathscr F}$ of an exhaustion of $\Sigma_\infty$, there exists a limiting global solution that converges to $F_\infty:\Sigma_\infty\to(M,g)$ by the Arzelà--Ascoli theorem in the setting of pointed Riemannian manifolds.

\begin{proposition}\label{prop-apl}
Let $(M,g)$ be an $(n\geq3)$-dimensional complete Riemannian manifold with bounded geometry, and let $\mathscr F :=\{\Sigma_t; \; t\in [0,T)\}$ be the MCF of $\Sigma$ in a gradient shrinking Ricci--Bourguignon soliton background which develops a singularity of type-I. Consider $\widetilde{\mathscr F}:=\{\widetilde \Sigma_s;\; s\in [-\log T,\infty)\}$ in $(M,g).$ Let $\{s_j\}_{j=1}^\infty$ be a divergent increasing sequence, and let $\{p_j\}_{j=1}^\infty$ be  a sequence of points in $\Sigma$ (assume Remark~\ref{unibou}, if $M$ is noncompact). Then, there exist subsequences $\{s_{j_k}\}$ and $\{p_{j_k}\}$ so that the family of immersions map $\widetilde{F}_{s_{j_k}} : \Sigma \to (M,g)$ from pointed manifolds $(\Sigma,\widetilde{F}_{s_j}^{*}g,p_{j_k})$ converges to an immersion map $F_\infty : \Sigma_\infty \to (M,g)$ from an $(n - 1)$-dimensional complete pointed Riemannian manifold $(\Sigma_\infty, F_\infty^*g, p_\infty)$ in the $C^\infty$ Cheeger--Gromov sense.
\end{proposition}

\begin{proof}
We will make this proof in a way that allows a generalization to the more general context of bounded geometry. Note that these hypotheses for the compact case are automatically satisfied. 

Take a sequence $\{\widetilde \Sigma_{s_j}\}_{j=1}^\infty$ in $\widetilde{\mathscr F}$ and points $\{p_j\}_{j=1}^\infty$ in $\Sigma$, and denote the Riemann curvature tensor of each $\widetilde \Sigma_{s_j}$ by $\widehat{\operatorname{Rm}}(\widetilde F_{s_j}^* g)$. 

Firstly, we will show that the sequence of pointed manifolds $\{(\Sigma,\widetilde F^*_{s_j}g,p_j)\}$ converges to some complete pointed Riemannian manifold $(\Sigma_{\infty}, h_\infty, p_\infty)$ in the $C^\infty$ Cheeger--Gromov sense up to subsequence. 

Since $(M,g)$ has bounded geometry, there exist positive constants $D_p$ and $\eta$ such that
\begin{equation}\label{Aux-GEq-WL1}
|\nabla^p \operatorname{Rm}_g|
\leq D_p \quad \mbox{and} \quad \operatorname{inj}(M,  g) \geq \eta > 0.
\end{equation}
for every integer $p\geq 0$. By Proposition~\ref{Underthesamesatupofproposition01} there are positive constants $C_p$ which does not depend on $s_j$ such that
\begin{equation}\label{Aux-GEq-WL2}
\big|\widehat\nabla_g^p \mathcal A(\widetilde F_{s_j})\big| \leq C_p.
\end{equation}
Thus, we are able to apply Thm.~2.1 by Chen and Yin \cite{chen2007uniqueness} which guarantees the existence of a positive constant $\delta = \delta(C_0, D_0, \eta, n)$ such that the injectivity radius of each $\big(\Sigma,\widetilde F^*_{s_j}g\big)$ satisfies
\begin{equation*}
\operatorname{inj}\big(\Sigma,\widetilde F^*_{s_j}g\big) \geq \delta > 0.
\end{equation*}
Moreover, from \eqref{Aux-GEq-WL1}, \eqref{Aux-GEq-WL2}, the Gauss equation and its iterated derivatives, we obtain positive constants $\widetilde C_p$ that also do not depend on $s_j$ such that
\begin{equation*}
\big|\widehat\nabla^p \widehat{\operatorname{Rm}}(\widetilde F_{s_j}^* g)\big| \leq \widetilde C_p
\end{equation*}
for every integer $p\geq0$. Hence, by Arzelà--Ascoli theorem, there exists a subsequence $\{(\Sigma, \widetilde  F_{s_{j_k}}^* g, p_{j_k})\}$ which converges to some complete pointed Riemannian manifold $(\Sigma_\infty, h_\infty, p_\infty)$, i.e., there exist an exhaustion $\{U_{j_k}\}_{k=1}^\infty$ of $\Sigma_{\infty}$ with $p_{\infty} \in U_{j_k}$ and diffeomorphisms $\Psi_{j_k}:U_{j_k}\to V_{j_k}:= \Psi_{j_k}(U_{j_k}) \subset \Sigma$ with $\Psi_{j_k}(p_\infty) = p_{j_k}$ such that $\Psi^*_{j_k}(\widetilde F^*_{s_{j_k}} g)$ converges in $C^{\infty}$ to $h_\infty$ uniformly on compact sets in $\Sigma_{\infty}$.

Secondly, we will prove the existence of an immersion map $F_{\infty} : \Sigma_\infty\to (M ,g)$. For this, let $\Theta: (M, g) \to (\mathbb R^d, g_{\rm st})$ be a Nash isometric embedding in some higher-dimensional Euclidean space such that, for every integer $j \geq 0$, the norm $|\nabla_g^{j} \mathcal A(\Theta)| \leq \overline D_{j}$, for some constants $\overline D_{j}>0$, where $\mathcal A(\Theta)$ is the second fundamental form of $\Theta$.
\begin{remark}\label{unibou}
At this point of the proof, we are assuming that the sequence $\{\widetilde F_{s_{j_k}}(p_{j_k})\}$ and the norms $|\nabla_g^{j} \mathcal A(\Theta)|$ are uniformly bounded in $M$, which are actually true in the compact case. Note, however, that we will need to prove or assume these facts for the noncompact case.
\end{remark}
Since we are working in the context of bounded geometry and each norm $|\nabla_g^{j} \mathcal A(\Theta)|$ is uniformly bounded, by setting 
\begin{equation*}
\overline F_{s_{j_k}} := \Theta \circ \widetilde F_{s_{j_k}} \circ \Psi_{j_k} :U_{j_k} \to (\mathbb R^d, g_{\rm st}),
\end{equation*}
we can prove that $\big|\nabla^p \overline F_{s_{j_k}}\big| \leq C_p$ for some positive constants $ C_p$ which does not depend on $s_{j_k}$, for every integer $p\geq0.$ The proof of this fact is by induction and follows the same steps as done in \cite[Appendix~C, pg.~42]{yamamoto2015meancurvature}. 

Then, by a standard argument as the Arzelà--Ascoli theorem, there exists a smooth map $\overline F_\infty: \Sigma_{\infty} \rightarrow (\mathbb R^d, g_{\rm st})$ such that the sequence of immersions $\overline F_{s_{j_k}} : U_{j_k} \to (\mathbb R^d,g_{\rm st})$ converges to $\overline F_{\infty} : \Sigma_{\infty}\rightarrow (\mathbb R^d,g_{\rm st})$ up to a subsequence. By definition of $C^{\infty}$ convergence, we get $h_\infty =\overline F^*_\infty g_{\rm st}$ on $\Sigma_\infty$. This implies that $\overline F_\infty: \Sigma_{\infty} \rightarrow \mathbb R^d$ is an isometric immersion map and $F_\infty := \Theta^{-1} \circ \overline F_{\infty}: \Sigma_\infty \to (M, g)$ is the required immersion of the theorem with $F_\infty^*g = (\Theta^{-1} \circ \overline F_{\infty})^*\Theta^*g_{\rm st} = h_\infty$.
\end{proof}

\begin{proof}[\bf Proof of Theorem~\ref{application:Huisken}]
We begin by taking a divergent increasing sequence $\left\{s_{j}\right\}_{j=1}^{\infty}$ and a sequence of points  $\{p_j\}_{j=1}^\infty$  in $\Sigma$. From Proposition~\ref{prop-apl}, there exist subsequences ${s_j}_k$ and ${p_j}_k$ so that the family of immersions map $\widetilde  F_{s_{j_{k}}}:\Sigma\to (M,g)$ from pointed manifolds $(\Sigma,p_{j_k})$ converges to an immersion map $F_{\infty}:\Sigma_{\infty}\to (M, g)$ from an $(n-1)$-dimensional compact pointed Riemannian manifold $(\Sigma_{\infty}, F_{\infty}^{*}g, p_\infty)$ in the $C^\infty$ Cheeger--Gromov sense. 

Next, we will prove that $\Sigma_{\infty}$ is a $f_{\infty}-$minimal hypersurface of $(M,g),$ where $f_\infty = f \circ F_\infty$. We denote $\widetilde{F}_{s_{j_k}}$ by $\widetilde {F}_{s_{k}}$ for short.  To simplify the notation we can take an exhaustion $\{U_k\}_{k=1}^\infty$ of $\Sigma_{\infty}$ with $p_{\infty} \in U_k$ and diffeomorphisms $\Psi_k : U_k \to V_k := \Psi_k(U_k) \subset \Sigma$ with $\Psi_k(p_\infty) = p_k$ such that $\Psi^*_k(\widetilde F^*_{s_{k}} g)$ converges in $C^{\infty}$ to $F^*_\infty  g$ uniformly on compact sets in $\Sigma_{\infty}$. Furthermore, the sequence of maps $\widetilde F_{s_{k}} \circ\Psi_k: U_k \rightarrow (M,g)$ converges in $C^{\infty}$ to $F_\infty: \Sigma_{\infty} \rightarrow(M,g)$ uniformly on compact sets in $\Sigma_{\infty}.$ So, for any compact set $K \subset \Sigma_\infty$ there exists $k_0$ such that $K \subset U_k$, for all $k \geq k_0$, and $\widetilde F_{s_{k}}\circ\Psi_k :U_k \rightarrow (M,g)$ converges to $F_\infty : \Sigma_\infty \to (M,g)$ in $C^\infty$ uniformly on $K$. Thus,
\begin{eqnarray*}
&&\int_K \left[H(\widetilde F_{s_{k}} \circ \Psi_k) + \nu_k\left(f\circ (\widetilde F_{s_{k}}\circ \Psi_k)\right) \right]^2 e^{- f\circ (\widetilde F_{s_{k}}\circ \Psi_k)} dA_{(\widetilde F_{s_{k}} \circ \Psi_k)^*g}\nonumber\\
&&\to \int_K \left( H( F_\infty) + \nu_\infty f_\infty\right)^ 2e^{- f_\infty }dA_{F^*_\infty g} 
\end{eqnarray*}
as $k\to \infty$, and 
\begin{eqnarray*}
&&\int_K \left[ H(\widetilde F_{s_{k}} \circ \Psi_k) +  \nu_k f (\widetilde F_{s_{k}}\circ\Psi_k)\right]^2e^{-f\circ (\widetilde F_{s_{k}}\circ \Psi_k)}dA_{(\widetilde F_{s_{k}} \circ\Psi_k)^* g}\\
&=&\int_{\Psi_k(K)}\left[H(\widetilde F_{s_{k}}) + \nu_k  (f\circ \widetilde F_{s_{k}}) \right]^2 e^{- f\circ \widetilde F_{s_{k}}} dA_{\widetilde F^*_{s_{k}}  g}\\
&\leq& \int_{\Sigma} \left[ H(\widetilde F_{s_{k}}) + \nu_k (f \circ \widetilde F_{s_{k}})\right]^2 e^{- f \circ\widetilde F_{s_{k}}} dA_{\widetilde F^*_{s_{k}}  g}.
\end{eqnarray*}
Hence, it is enough to prove the following: 
\begin{equation}\label{lastconv3}
\int_\Sigma [H(\widetilde F_{s_{k}})+ \nu_k( f \circ \widetilde F_{s_{k}})]^2 e^{-f\circ\widetilde F_{s_{k}}} dA_{\widetilde F_{s_{k}}^*g}\to 0
\end{equation}
as $k\to \infty$. We will argue by contradiction. Assume that there exist a constant $\delta>0$ and a subsequence $\{\ell\}\subset \{k\}$ with $\ell\to\infty$ such that 
\begin{equation*}
\int_\Sigma [H(\widetilde F_{s_{\ell}})+\nu_\ell(f \circ \widetilde F_{s_{\ell}})]^2 e^{-f\circ\widetilde F_{s_{\ell}}} dA_{\widetilde F_{s_{\ell}}^*g}\geq\delta.
\end{equation*}
Then
\begin{equation*}
\int_\Sigma [H(\widetilde F_s)+ \nu_s(f \circ \widetilde F_s)]^2 e^{-f\circ\widetilde F_s} dA_{\widetilde F_s^*g} \geq \frac{\delta}{2}, 
\end{equation*}
for $s\in[s_{\ell},s_{\ell}+\frac{\delta}{2C'}]$, 
where we used Lemma~\ref{Lem:4.10} and $C'$ is the constant appeared in that lemma. Hence, 
\begin{equation*}
\int_{-\log T}^\infty \int_\Sigma [H(\widetilde F_s)+ \nu_s(f \circ \widetilde F_s)]^2 e^{-f\circ\widetilde F_s} dA_{\widetilde F_s^*g} ds =\infty. 
\end{equation*}
On the other hand, by the monotonicity formula in Lemma~\ref{Prop:4.6}
\begin{equation*}
\frac{d}{ds} \int_\Sigma e^{-f\circ\widetilde F_s} dA_{\widetilde F_s^*g}=-\int_\Sigma [H(\widetilde F_s)+ \nu_s(f \circ \widetilde F_s)]^2 e^{-f\circ\widetilde F_s} dA_{\widetilde F_s^*g}\leq 0. 
\end{equation*}
Thus, the weighted volume 
\begin{equation*}
\int_\Sigma e^{-f\circ\widetilde F_s} dA_{\widetilde F_s^*g} 
\end{equation*}
is monotone decreasing and nonnegative. Therefore, it converges to some value 
\begin{equation*}
\alpha:=\lim_{s\to\infty}\int_\Sigma e^{-f\circ\widetilde F_s} dA_{\widetilde F_s^*g}<\infty, 
\end{equation*}
and then we obtain the following contradiction:
\begin{equation*}
\int_{-\log T}^\infty \int_\Sigma [H(\widetilde F_s)+ \nu_s(f \circ \widetilde F_s)]^2 e^{-f\circ\widetilde F_s} dA_{\widetilde F_s^*g}
=-\alpha +\int_\Sigma e^{-f\circ\widetilde F_a} dA_{\widetilde F_a^*g}<\infty,  
\end{equation*}
where $a:=-\log T$, which proves \eqref{lastconv3}. With this, the proof of the theorem is complete.
\end{proof}

\section{\bf The noncompact case}\label{Sec-ncompcase}
In this section, we address the case of complete noncompact Riemannian manifolds $(M,g)$ with some additional uniformity conditions. We begin with a brief discussion on the reduced distance along the Ricci--Bourguignon flow, initially defined by Perelman in the Ricci flow setting (see~\cite[Sect.~7]{perelman2002entropy}).  

Let $\overline g(t)$ be a shrinking self-similar solution of the Ricci--Bourguignon flow in $M \times [0,T)$.  For any smooth curve $\gamma:[t_{1},t_{2}]\to M$ with $0\leq t_{1}<t_{2}<T$, consider the $\mathcal{L}$-length of $\gamma$ by 
\begin{equation*}
\mathcal{L}(\gamma):=\int_{t_{1}}^{t_{2}}\sqrt{t_{2}-t}\bigl( (1-n\rho)R_{\overline g}+|\dot{\gamma}|^2 \bigr)dt, 
\end{equation*}
where $|\dot{\gamma}|$ is the norm of $\dot{\gamma}(t)$ measured by $\overline g$. 
For a fixed point $(q_{2},t_{2})$ in the space-time $M\times [0,T)$, we get the reduced distance  
\begin{equation*}
\ell_{q_{2},t_{2}}:M\times [0,t_{2})\to\mathbb{R}
\end{equation*}
based at $(q_{2},t_{2})$ by
\begin{equation*}
\ell_{q_{2},t_{2}}(q_{1},t_{1}):=\frac{1}{2\sqrt{t_{2}-t_{1}}}\inf_{\gamma}\mathcal{L}(\gamma), 
\end{equation*}
where the infimum is taken over all smooth curve $\gamma:[t_{1},t_{2}]\to M$ with $\gamma(t_{1})=q_{1}$ and $\gamma(t_{2})=q_{2}$. 

In what follows, we assume that there exists a Nash isometric embedding $\Theta: (M, g) \to (\mathbb R^d, g_{\rm st})$ in some higher-dimensional Euclidean space such that for every integer $j \geq 0,$ the second fundamental form $\mathcal A(\Theta)$ of $\Theta$ satisfies
\begin{equation}\label{CAdd-Thm3}
|\nabla_g^{j} \mathcal A(\Theta)| \leq \overline D_{j}
\end{equation}
for some constants $\overline D_{j}>0$. We observe that, under this assumption, $(M,g)$ must have bounded geometry by the Gauss equation (and its iterated derivatives) and Thm.~2.1 of \cite{chen2007uniqueness}.

\begin{theorem}\label{Thmnoncommpact}
Let $(M,g)$ be an $(n\geq3)$-dimensional complete noncompact Riemannian manifold
satisfying \eqref{CAdd-Thm3}, and let $\mathscr F:=\{\Sigma_t; \; t\in [0,T)\}$ be the MCF of $\Sigma$ in a gradient shrinking Schouten soliton background which develops a singularity of type-I. Consider $\widetilde{\mathscr F}:=\{\widetilde \Sigma_s;\; s\in [-\log T,\infty)\}$ in $(M,g)$ such that $\phi\circ \widetilde{F}_{s}$
remains uniformly bounded on $[-\log T,\infty).$  In addition, assume that there exists a point $q_{0}\in \Sigma$ such that the reduced distance $\ell_{F_{t}(q_{0}),t}$ converges pointwise to $f$ (as $t\to T$) on $M\times[0,T).$ Then, for any divergent increasing sequence $\left\{s_{j}\right\}_{j=1}^{\infty},$ there exists a subsequence $\{s_{j_{k}}\}$ such that the family of immersion maps $\widetilde F_{s_{j_k}}: \Sigma\to (M,g)$ from pointed manifold $\big(\Sigma ,q_0\big)$ converges to an immersion map $F_\infty: \Sigma_\infty \to (M,g)$ from an $(n-1)$-dimensional complete pointed Riemannian manifold $\left(\Sigma_{\infty}, F_\infty ^* g, q_{\infty}\right)$ in the $C^\infty$ Cheeger--Gromov sense. Furthermore, $(\Sigma_{\infty}, F_\infty^* g)$ is an $f_\infty-$minimal hypersurface of  $(M, g),$ where $f_\infty = f \circ F_\infty$. 
\end{theorem}

\begin{proof}
As we had already mentioned in the proof of Proposition~\ref{prop-apl}, it is enough to show that $\{\widetilde F_{s_{j}}(q_{0})\}_{j = 1}^\infty$ is a bounded sequence in $(M,g)$, since the remainder of the proof is the same as in the compact case. 

We begin by taking $t_1$, $t_2$ with $0\leq t_{1}<t_{2}<T$, and $\{F_t(q_0)\}_{t\in[t_1, t_2]}$ as a curve joining $F_{t_1}(q_0)$ and $F_{t_2}(q_0)$. Thus, 
\begin{eqnarray*}
\ell_{F_{t_2}(q_0),t_2}(F_{t_1}(q_0),t_1)&\leq& \frac{1}{2\sqrt{t_2-t_1}}\int_{t_{1}}^{t_{2}}\sqrt{t_{2}-t}\left( (1-n\rho)R_{\overline g} + \left|\frac{\partial F_{t}}{\partial t}\right|^2 \right)dt\\
&=& \frac{1}{2\sqrt{t_{2}-t_{1}}}\int_{t_{1}}^{t_{2}}\sqrt{t_{2}-t}\bigl( (1-n\rho)R_{\overline g} + H_{\overline g}^2 \bigr)dt.
\end{eqnarray*}
Since $\mathscr F$ develops a singularity of type-I, we have $(T-t)H_{\overline g}^2$ is bounded. Moreover, bounded geometry assumption and the fact that $(T-t)R_{\overline g}=R_g$ imply $(1-n\rho)R_{\overline g} + H_{\overline g}^2\leq \frac{C}{T-t}$ for some positive constant $C$ which does not depend on $t$. Hence,
\begin{eqnarray*}
&&\ell_{F_{t_2}(q_{0}),t_{2}}(F_{t_{1}}(q_0),t_{1})\\
&\leq&  \frac{C}{2\sqrt{t_{2}-t_{1}}}\int_{t_{1}}^{t_{2}}\frac{\sqrt{t_{2}-t}}{T-t}dt
\leq   \frac{C}{2\sqrt{t_{2}-t_{1}}}\int_{t_{1}}^{t_{2}}\frac{1}{\sqrt{T-t}}dt
\leq  C\frac{\sqrt{T-t_{1}}}{\sqrt{t_{2}-t_{1}}}. 
\end{eqnarray*}
By assumption $\ell_{F_t(q_0),t}$ converges pointwise to $f$ (as $t\to T$) on $M\times[0,T)$ and by taking the limit as $t_{2}\to T$, we have 
$f(F_{t_1}(q_0),t_1)\leq C$. As $f(F_t(q_0),t)=\overline f(F_t(q_0))=f(\widetilde F_s(q_0))$, one has $f(\widetilde F_s(q_0))\leq C$
for all $s\in[-\log T,\infty)$. It follows from Borges~\cite[Thm.~1.2]{borges2022on} that there exist positive constants $A_1$ and $A_2$ depending on $\lambda$ and the geometry of the soliton on the unity ball $B_{q_0}(1)$ such that
\begin{equation*}
\frac{\lambda}{4}\left( r - A_1\right)^2 + f_0 \leq f \leq \lambda \left(r+A_2\right)^2 + f_0
\end{equation*}
on $M$, where $r(q)=d_g(q_0, q)$ is the distance function from any fixed point $q_0\in M$ with $r(q)>2$ and $f_0=\min_{p \in M} f(p).$ Therefore
\begin{equation*}
d_g(q_0, \widetilde F_s(q_0))\leq \sqrt{\frac{4}{\lambda}\left(C - f_0\right)}+A_1, 
\end{equation*}
from which that $\{\widetilde F_{s_{j}}(q_{0})\}_{j = 1}^\infty$ is bounded in $(M,g)$ and the theorem follows.
\end{proof}

\section{\bf Construction of mean curvature solitons in a gradient Ricci--Bourguignon soliton Background}\label{example:section} 

In this section, we show how to construct explicit examples of MCF in a gradient Ricci--Bourguignon soliton background. For explicit examples of MCF in a Ricci flow background, see the work by Yamamoto~\cite{yamamoto2018examples}.

Let $g$ be a complete Riemannian metric on an $(n\geq3)$-dimensional smooth manifold $M$ satisfying
\begin{align}\label{Bizu1}
\operatorname{Ric}_g + \operatorname{Hess}_g f = (\lambda + \rho R_g)g,
\end{align}
for some real constants $\lambda$ and $\rho$, and a nonconstant smooth potential function $f$ on $M$.

A direct calculation shows that a gradient Ricci--Bourguignon soliton, i.e., a self-similar solution $g(t)=\sigma(t)\psi_t^*g$ of \eqref{eq-RB}, is characterized by \eqref{Bizu1}. Here, $\psi_t$ is the flow of the vector field $\nabla_g f/\sigma(t)$ with parameter $t$ chosen such that $\sigma(t):=1-2\lambda t$ is a positive constant.  Two fundamental references for the study of gradient Ricci--Bourguignon solitons are the works of Catino, Mazzieri and Mongodi~\cite{catino2015rigidity} and Catino and Mazzieri~\cite{catino2016gradient}.

Let $g_{ij} = \varphi^{-2}\delta_{ij}$ be a Riemannian metric on $\mathbb{R}^n,$ where $\delta_{ij}$ stands for the Euclidean metric and $\varphi$ is a nonzero smooth function on $\mathbb{R}^n.$ It is known that 
\begin{align*}
(\operatorname{Ric})_{ij} =\frac{1}{\varphi^2}\Big((n - 2)\varphi \varphi_{x_i x_j} + \Big(\varphi\sum_k \varphi_{x_kx_k} - (n - 1)\sum_k \varphi_{x_k}^2\Big)\delta_{ij}\Big),
\end{align*}
where the two terms in the second part of this equation are evaluated using the metric $\delta_{ij},$ see \cite[p. 59]{besse2007einstein}. Therefore,
\begin{align}\label{ScalarConform}
R_g = 2(n-1)\varphi\sum_k \varphi_{x_kx_k} - n(n-1)\sum_k \varphi_{x_k}^2,
\end{align}
and the following equations are valid
\begin{align*}
(\operatorname{Hess}_g h)_{ij} &= h_{x_ix_j} +\frac{\varphi_{x_j}}{\varphi} h_{x_i} +\frac{\varphi_{x_i}}{\varphi}
h_{x_j} \quad \rm{for} \quad i \neq j,\\
(\operatorname{Hess}_g h)_{ii} &= h_{x_ix_i} + 2\frac{\varphi_{x_i}}{\varphi} h_{x_i} -\sum_k\frac{\varphi_{x_k}}{\varphi} h_{x_k} \quad \rm{for}\quad i = j,  
\end{align*}
for any smooth function $h$ on $\mathbb{R}^n$. 

We seek solutions of \eqref{Bizu1} of the form $f(\xi)$, depending only on $\xi=\sum_{i=1}^n\alpha_i x_i$, where $\alpha=(\alpha_1,\ldots,\alpha_n)\in\mathbb{R}^n$ is nonzero. Since $\sum_{i=1}^n\alpha_i^2\neq0$, we may assume, without loss of generality, that $\sum_{i=1}^n\alpha_i^2=1$. The following proposition characterizes the system of ordinary differential equations satisfied by such solutions and provides the parameters to construct gradient Ricci--Bourguignon solitons on $\mathbb{R}^n$. 

\begin{proposition}\label{Feitosa-FreitasFilho-Gomes}
Let $\mathbb{R}^n$ be a Euclidean space with  metric $g_{ij} = \varphi^{-2}(\xi)\delta_{ij}$, where $\varphi (\xi )$ is a nonzero smooth smooth function on $\mathbb{R}^n$ and $\xi  = \sum_{i = 1}^n \alpha_ix_i$ with $\sum_{i = 1}^n \alpha_i^2 = 1$. We can obtain a smooth function  $f (\xi )$ satisfying~\eqref{Bizu1} by solving the following equations
\begin{align}\label{rn1}
(n - 2) \dfrac{\varphi''}{\varphi} + f''+2\dfrac{\varphi'}{\varphi} f'=0
\end{align}
and
\begin{align}\label{rn2}
(1-2(n-1)\rho)\dfrac{\varphi''}{\varphi} + (n - 1)(n\rho -1) \Big(\dfrac{\varphi'}{\varphi}\Big)^2 -\dfrac{\varphi'}{\varphi}f'= \dfrac{\lambda}{\varphi^2}.   
\end{align}
Here, the superscript $'$ denotes differentiation with respect to $\xi$.
\end{proposition}

\begin{proof}
We need to analyze \eqref{Bizu1} in two cases. For $i \neq j$, it rewrites as
\begin{align}\label{Holds}
(n - 2) \frac{\varphi_{x_ix_j}}{\varphi} + f_{x_i x_j} +\frac{\varphi_{x_j}}{\varphi} f_{x_i} +\frac{\varphi_{x_i}}{\varphi}f_{x_j}  = 0
\end{align}
and, for all $i$, as
\begin{eqnarray}\label{BIZU1}
\frac{\lambda}{\varphi^2}&=& (n - 2) \frac{\varphi_{x_ix_i}}{\varphi} + f_{x_i x_i} +2\frac{\varphi_{x_i}}{\varphi} f_{x_i} \\  
\nonumber&&+ \!\sum_k\!\Big( (1-2(n-1)\rho)\frac{ \varphi_{x_k x_k} }{\varphi}
+\! (n - 1)(n\rho -1)\frac{\varphi^2_{x_k}}{\varphi^2} -\!\frac{\varphi_{x_k}}{\varphi}f_{x_k}\Big).
\end{eqnarray}
We now assume that the functions $\varphi(\xi)$ and $f(\xi )$ depend only on $\xi.$ Hence, we have $\varphi_{x_i}= \varphi'\alpha_i$ and $\varphi_{x_ix_j} = \varphi''\alpha_i\alpha_j.$ Using the same reasoning for $f$, equations~\eqref{Holds} and~\eqref{BIZU1} become, for  $i \neq j$,
\begin{align}\label{Sys:Bizu1}
\left[(n - 2)\frac{\varphi''}{\varphi} + f'' +2\frac{\varphi'}{\varphi} f' \right]\alpha_i\alpha_j= 0
\end{align}
and, for all $i,$
\begin{align}\label{Sys:Bizu2}
 \frac{\lambda}{\varphi^2}
=& \Big[(n - 2) \frac{\varphi''}{\varphi} + f''+2\frac{\varphi'}{\varphi} f'\Big]{\alpha_i}^2\\  
&+ \!\sum_k\!\Big( (1-2(n-1)\rho)\frac{ \varphi'' }{\varphi}{\alpha_k}^2
+\! (n - 1)(n\rho -1)\frac{{\varphi'}^2}{\varphi^2}{\alpha_k}^2 -\!\frac{\varphi'}{\varphi}f'{\alpha_k}^2\Big).\nonumber
\end{align}
Since $\sum_{k = 1}^n \alpha_k^2 = 1$, equations \eqref{rn1} and \eqref{rn2} follow from \eqref{Sys:Bizu1} and \eqref{Sys:Bizu2}.
\end{proof}
\begin{remark}
For constructing a family of mean curvature solitons $\mathscr G$ in the Ricci--Bourguignon flow background on $\mathbb R^n$ with Riemannian metric $\varphi^{-2}(\xi)\delta_{ij}$, it is enough to find parameters that satisfy~\eqref{rn1} and~\eqref{rn2} and an $f$-minimal hypersurface $\Sigma$ in this geometric ambient space. Indeed, it follows immediately from the characterization of a gradient soliton for the Ricci--Bourguignon flow, Theorem~\ref{NazasCharacterization} and Proposition~\ref{Feitosa-FreitasFilho-Gomes}. It is also important to note that an $f$-minimal hypersurface is characterized as a critical point of the $f$-weighted area functional.
\end{remark}

In what follows, we work in the geometric setting of Proposition~\ref{Feitosa-FreitasFilho-Gomes}, which will be instrumental in our arguments. We take $\mathbb R^n = \mathbb R \times \mathbb R^{n-1}$ with the Riemannian metric $g$ given by $g_{ij} = \varphi^{-2}(\xi)\delta_{ij}$, where $\xi  = \sum_{i = 1}^n \alpha_ix_i$ and $(\alpha_1,\ldots,\alpha_n)$ is a unit vector in $\mathbb{R}^n$. Consider $\mathbb T^{n-1}=\mathbb R^{n-1}/\mathbb Z^{n-1}$, and let $\pi:\mathbb R^n\to \mathbb R\times\mathbb T^{n-1}$ be the canonical projection map. Take the quotient metric induced by $\pi$ to obtain a local isometry between  $\big(\mathbb R^n,g\big)$ and $(\mathbb R\times \mathbb T^{n-1},g_0)$, where $g_0=\varphi^{-2}(\xi) \left(d\xi^2+g_{\mathbb T^{n-1}}\right)$ and $g_{\mathbb T^{n-1}}$ is the flat quotient standard metric, provided that $\varphi(\xi)$ is invariant under the action of $\mathbb Z^{n-1}.$  

\begin{example}\label{Ex-Agila-Gomes}
The parameters $\varphi(x)=\sech x,$ $\lambda=\rho=1/3$, and $f(x)= \frac{2}{3}\ln(\cosh x) +\frac16\cosh^2x$ define $g_0=\sech^{-2}x(dx^2+g_{\mathbb T^2})$ as a complete gradient shrinking traceless Ricci--Bourguignon soliton on $\mathbb{R}\times\mathbb{T}^2$ with scalar curvature $R_{g_0}=-2\sech^2 x(1+\sech^2x)$, satisfying $-4\leq R_{g_0}<0$. Moreover, $\Sigma_0=\{0\}\times\mathbb T^2$ is a compact $f$-minimal hypersurface of $\mathbb{R}\times\mathbb{T}^2$.
\end{example}

Indeed, the first part follows from Proposition~\ref{Feitosa-FreitasFilho-Gomes}, and $R_{g_0}$ can be computed directly by \eqref{ScalarConform}. Moreover, the divergent curve criterion immediately implies that $g_0$ is geodesically complete. For the remainder, note that each slice $\Sigma_x=\{x\}\times\mathbb{T}^2$ has induced metric
$g_{\Sigma_x}=\varphi^{-2}(x)g_{\mathbb T^2}$, and then $\operatorname{Area}_f(\Sigma_x) = e^{-f(x)}\varphi(x)^{-2}\operatorname{Area}(\mathbb T^2).$ Thus, 
\begin{equation}
\frac{d}{dx}\operatorname{Area}_f(\Sigma_x) = e^{-f(x)}\varphi(x)^{-2}\big(f'(x)+2\varphi'(x)/\varphi(x)\big)\operatorname{Area}(\mathbb T^2).
\end{equation}
So, by the variational characterization of $f$-minimal hypersurfaces, $\Sigma_x$ is $f$-minimal if and only if $f'(x)+2\varphi'(x)/\varphi(x)=0$. In particular, at $x=0$, we have $f'(0)+2\varphi'(0)/\varphi(0)=0$ which is equivalent to $\Sigma_0$ being $f$-minimal.

\begin{remark}
The approach used in Example~\ref{Ex-Agila-Gomes} differs from that used in Example~2 of Agila and the first author~\cite{agilagomes}, but leads to the same result.
\end{remark}

\section{\bf Acknowledgments}
José N.V. Gomes has been partially supported by the Conselho Nacional de Desenvolvimento Científico e Tecnológico (CNPq), Grant 309825/2026-1, and by the Fundação de Amparo à Pesquisa do Estado de São Paulo (FAPESP), Grants 2024/00923-6, 2023/11126-7, and 2022/16097-2. Matheus Hudson has also been partially supported by FAPESP, Grants 2023/13921-9 and 2024/20015-7. Raquel Lemos acknowledges support from CNPq, Grant 153758/2022-8.

\end{document}